\documentclass[12pt, eqno]{article}

\usepackage{amsmath,amssymb,latexsym}

\newcommand{\F}{{\mathbb F}}
\newcommand{\N}{{\mathbb N}}
\newcommand{\C}{{\mathbb C}}
\newcommand{\R}{{\mathbb R}}
\newcommand{\Q}{{\mathbb Q}}
\newcommand{\Z}{{\mathbb Z}}
\newcommand{\CP}{{\mathbb CP}}
\newcommand{\RP}{{\mathbb RP}}

\newtheorem{theorem}{Theorem}[section]

\newtheorem{corollary}[theorem]{Corollary}

\newtheorem{example}[theorem]{Example}
\newtheorem{remark}[theorem]{Remark}

\newtheorem{lemma}[theorem]{Lemma}

\newfont{\Bb}{msbm10 scaled\magstephalf}

\begin{document}
\title{On the degenerated Arnold-Givental conjecture}
\author{Guangcun Lu\thanks{The first author is partially sponsored by the NNSF 10671017
and 11271044  of China and the Program for New Century Excellent Talents of the Education
Ministry of China.}\\
Department of Mathematics,\\
Beijing Normal University, China\\  \vspace{2mm}
(gclu@bnu.edu.cn)\\}

\date{First version, May 30, 2008\\
Revised version, August 6, 2018}
\maketitle \vspace{-0.1in}

\abstract{We present another view dealing with
the  Arnold-Givental conjecture on a real symplectic manifold
 $(M, \omega, \tau)$  with  nonempty and compact real part $L={\rm Fix}(\tau)$.
 For given $\Lambda\in (0, +\infty]$ and $m\in\N\cup\{0\}$
  we show the equivalence of the following two claims: (i)
 $\sharp(L\cap\phi^H_1(L))\ge m$ for any Hamiltonian function $H\in C_0^\infty([0, 1]\times M)$ with
Hofer's norm $\|H\|<\Lambda$; (ii)  $\sharp {\cal P}(H,\tau)\ge m$
for every $H\in C^\infty_0(\R/\Z\times M)$ satisfying
$H(t,x)=H(-t,\tau(x))\;\forall (t,x)\in\mathbb{R}\times M$ and
with Hofer's norm $\|H\|<2\Lambda$, where
${\cal P}(H, \tau)$ is the set of all $1$-periodic solutions of $\dot{x}(t)=
X_{H}(t,x(t))$ satisfying $x(-t)=\tau(x(t))\;\forall t\in\R$ (which are also called  brake orbits
sometimes). Suppose that $(M, \omega)$ is geometrical bounded for some $J\in{\cal J}(M,
\omega)$ with $\tau^\ast J=-J$ and has a rationality index $r_\omega>0$ or $r_\omega=+\infty$.
Using Hofer's method we  prove that if the Hamiltonian $H$ in (ii) above has  Hofer's norm
$\|H\|<r_\omega$ then  $\sharp(L\cap\phi^H_1(L))\ge\sharp {\cal P}_0(H,\tau)\ge
{\rm Cuplength}_{\F}(L)$ for $\F=\Z_2$, and further for $\F=\Z$ if $L$ is
orientable, where ${\cal P}_0(H,\tau)$ consists of all
contractible solutions in ${\cal P}(H,\tau)$. }
 \vspace{-0.1in}
\medskip

%%%%%%%%%%%%%%%%%%%%%%%%%%%%%%%%%%%%%%%%%%%%%%%%%%%%%%%%%%%%%%%%%%%%%%%%%%%%%%%%%%%%
%% The conclusion in $r_\omega=+\infty$ may be considered as a partial
%% affirmative answer to the  degenerated Arnold-Givental conjecture.
%%%%%%%%%%%%%%%%%%%%%%%%%%%%%%%%%%%%%%%%%%%%%%%%%%%%%%%%%%%%%%%%%%%%%%%%%%%%%%%%%%%%%%%%%%%

\section{Introduction}\label{sec:1}

A {\bf real symplectic manifold} is a triple $(M,\omega,\tau)$
consisting of a symplectic manifold $(M,\omega)$ and an
anti-symplectic involution $\tau$ on $(M,\omega)$, i.e.
$\tau^\ast\omega=-\omega$ and $\tau^2=id_M$.
The Marsden-Weinstein quotients of real Hamiltonian systems
provide a great deal of examples of such manifolds.
 Let ${\cal J}(M,\omega)$ denote the space of all $\omega$-compatible smooth
almost complex structures on $M$, and
$$
\R{\cal J}(M,\omega)=\{J\in {\cal J}(M,\omega)\,|\, J\circ
d\tau=-d\tau\circ J\},
$$
that is, $J\in\R{\cal J}(M,\omega)$ if and only if $\tau$ is
anti-holomorphic with respect to $J$. With the standard trick of
S\'evennec (see \cite[p.64]{McSa1}) one can prove that $\R{\cal
J}(M,\omega)$ is a separable Frech\'et submanifold of ${\cal
J}(M,\omega)$ which is nonempty and contractible (cf. \cite[Prop.
1.1]{Wel}). The fixed point set $L:={\rm Fix}(\tau)$ of $\tau$ is
called the {\bf real part} of $M$. Since $\tau$ is an isometry of
the natural Riemann  metric $g_J=\omega\circ(id_M\times J)$ for any
$J\in\R{\cal J}(M,\omega)$, $L$ is either empty or a Lagrange
submanifold (cf. \cite[p.4]{Vi}).\\

\noindent{\bf Arnold-Givental conjecture}(\cite{Gi}): Let
$(M,\omega,\tau)$ be a real  symplectic manifold of dimension $2n$,
and $L={\rm Fix}(\tau)$ be a nonempty compact submanifold without
boundary. Then for every Hamiltonian diffeomorphism $\phi$ on
$(M,\omega)$, it holds that
\begin{equation}\label{e:1.7}
\sharp\bigl(L\cap\phi(L)\bigr)\ge
\sum^n_{k=0}b_k(L,\Z_2)\quad\hbox{or}\quad\sum^n_{k=0}b_k(L,\Z)
\end{equation}
provided that $L$ and $\phi(L)$ intersect transversally, and that
\begin{equation}\label{e:1.8}
\sharp\bigl(L\cap\phi(L)\bigr)\ge {\rm
Cuplength}_{\Z_2}(L)\quad\hbox{or}\quad {\rm Cuplength}_{\Z}(L)
\end{equation}
generally. Hereafter the $\F$-cuplength of a paracompact topological space $X$
over an integral domain $\F$, ${\rm Cuplength}_{\F}(X)$, is defined
the supremum of natural numbers $k$ such that there exist cohomology
classes $\alpha_1,\cdots, \alpha_{k-1}$ in $H^\ast (X,\F)$ of
positive degree satisfying $\alpha_1\cup\cdots\cup\alpha_{k-1}\ne
0$.\vspace{2mm}

This conjecture is a special case of Arnold's conjecture on
Lagrangian intersections (\cite{Ar1, Ar2}).
 If $M$ is closed, the estimate (\ref{e:1.7}) in $\Z_2$-coefficients
follows from Floer \cite{Fl1} if $\pi_2(M,L)=0$,
\cite{Oh} if $L$ is a real form of compact Hermitian spaces with some assumptions on
the Maslov index, \cite{Laz} if $L$ is the strongly negative monotone, and \cite[Theorem H]{FuOOO} if $L$ is the semipositive,  and \cite{Fr} if $L$ is in Marsden-Weinstein
quotients. The estimate (\ref{e:1.8}) in
$\Z_2$-coefficients follows from Floer and Hofer \cite{Fl2, Ho2},
and Liu \cite{Liu}
if $(M, \omega)$ has positive rationality index $r_\omega$
and $\phi$ may be generated by $H\in
C^\infty([0, 1]\times M)$ with Hofer's norm $\|H\|<r_\omega/2$.
The estimates in (\ref{e:1.7}) and (\ref{e:1.8}) were obtained for
$(M, L)=(\CP^n, \RP^n)$ \cite{ChJi, Gi}. (The author \cite{Lu2} also
generalized the arguments in \cite{ChJi} to the case of weighted
complex projective spaces, which are symplectic orbifolds).

Arnold-Givental conjecture
contains  Arnold conjecture for
the symplectic fixed points (\cite{Ar1, Ar2}), which stated that for every Hamiltonian
diffeomorphism $\phi$ on a closed symplectic manifold $(M,\omega)$
the following estimates hold true,
\begin{eqnarray}
&&\sharp{\rm Fix}(\phi)\ge{\rm Cuplength}_{\F}(M),\label{e:1.9}\\
&&\sharp{\rm Fix}(\phi)\ge \sum^{\dim M}_{k=0}b_k(M;\F)\label{e:1.10}
\end{eqnarray}
if each $x\in {\rm Fix}(\phi)$ is nondegenerate in the sense that
 the tangent map $d\phi(x):T_xM\to T_xM$ has no eigenvalue $1$.
 After Floer \cite{Fl3} first invented Floer
homologies to prove  (\ref{e:1.10}) for monotone $(M,\omega)$ and $\F=\Z$, Fukaya-Ono \cite{FuO} and Liu-Tian
\cite{LiuT}  developed Floer homologies to affirm it for any closed symplectic
manifold $(M,\omega)$ and $\F=\Q$.

In this paper we consider a smooth time dependent Hamiltonian function $H:\R\times
M\to\R, (t,x)\mapsto H(t,x)=H_t(x)$ satisfying
\begin{equation}\label{e:1.1}
H_t(x)=H_{t+1}(x)\quad{\rm and}\quad H(t, x)=H(-t, \tau(x))\;\forall
(t,x)\in\R\times M.
\end{equation}
Such a Hamiltonian function $H$ is said to be {\bf $1$-periodic in
time and symmetric}. Let $X_{H_t}$ be defined by
$\omega(X_{H_t},\cdot)=-dH_t(\cdot)$. %the same as in [HoZe], converse as [McSa]%
Then $X_{H_t}=X_{H_{t+1}}$ and
\begin{equation}\label{e:1.2}
X_{H_{-t}}(x)=-d\tau(\tau(x))X_{H_t}(\tau(x))\;\;\forall
(t,x)\in\R\times M.
\end{equation}
 For $x_0\in M$ let $x:\R\to M$
be the solution  of
\begin{equation}\label{e:1.3}
\dot x(t)=X_{H_{t}}(x(t))
\end{equation}
through $x_0$ at $t=0$. Then both $y(t):=x(-t)$ and
$z(t):=\tau(x(t))$ are solutions of
$$\dot x(t)=d\tau(\tau(x(t))X_{H_{t}}(\tau(x(t))).$$
 So
$y=z$ if and only if $x_0=y(0)=z(0)=\tau(x(0))=\tau(x_0)$. We are
interested in
 those $1$-periodic solutions $x$ of the equation (\ref{e:1.3}) which satisfy
\begin{equation}\label{e:1.4}
 x(-t)=\tau(x(t))\;\forall t\in\R.
\end{equation}
Clearly, such a solution $x$ satisfies: $x(0), x(1/2)\in L$ and $x(1-t)=\tau(x(t))\;\forall t$.
A loop  $x\colon S^{1}=\R/\Z\to M$ satisfying (\ref{e:1.4}) is
called a {\bf $\tau$-reversible}. ($\tau$-reversible $1$-periodic solutions
are also called {\bf brake orbits} in literature.)
 Denote by
$$
{\cal P}(H, \tau)\;({\rm resp}.\; {\cal P}_0(H,\tau)\,)
$$
the set of all $\tau$-reversible $1$-periodic solutions (resp.
contractible $\tau$-reversible $1$-periodic solutions) of
(\ref{e:1.3}). Clearly, ${\cal P}(H, \tau)$ must be empty if $L=\emptyset$.
 Let $\phi^H_t:M\to M$ be the Hamiltonian  diffeomorphisms defined by
$$
\frac{d}{dt}\phi^H_t=X_{H_t}\circ\phi^H_t,\quad \phi^H_0=id_M.
$$
From (\ref{e:1.2}) it easily follows that
$\phi^H_t\circ\tau=\tau\circ\phi^H_{-t}\;\forall t\in\R$. Moreover,
it always holds that $\phi^H_{t+1}=\phi^H_t\circ\phi^H_1\;\forall
t\in\R$. So we get that
\begin{equation}\label{e:1.5}
\phi^H_1\circ\tau=\tau\circ(\phi^H_1)^{-1}.
\end{equation}
One also easily checks that the elements of ${\cal P}(H,\tau)$ are
one-to-one correspondence with points in $L\cap{\rm Fix}(\phi^H_1)$.
 So we have
\begin{equation}\label{e:1.6}
\sharp(L\cap{\rm Fix}(\phi^H_1))=\sharp{\cal
P}(H,\tau)\ge\sharp{\cal P}_0(H,\tau).
\end{equation}
Recall  Hofer's norm of a Hamiltonian function $H\in
C^\infty_0([0,1]\times M)$ is defined by
$$
\|H\|=\int^1_0[\sup_xH_t(x)-\inf_xH_t(x)]dt.
$$
Our first result is

\begin{theorem}\label{th:1.1}
 Let $(M,\omega,\tau)$ be a real  symplectic manifold of dimension $2n$, and the fixed point set $L={\rm
Fix}(\tau)$ be nonempty. Let $\Lambda\in (0, +\infty]$ and
$m\in\N\cup\{0\}$. Then the
following two claims are equivalent. \\
{\rm (i)} Every Hamiltonian diffeomorphism $\phi$ on $M$ generated
by a Hamiltonian function $H\in C_0^\infty([0, 1]\times M)$ with
$\|H\|<\Lambda$, satisfies
$$
\sharp(L\cap\phi(L))\ge m.
$$
{\rm (ii)} Every 1-periodic in time and symmetric $H\in
C^\infty_0(\R/\Z\times M)$ with $\|H\|<2\Lambda$,
satisfies
$$
\sharp {\cal P}(H,\tau)\ge m.
$$
\end{theorem}

\begin{remark}\label{rem:1.2}
{\rm The proof of ``(i)$\Longrightarrow$(ii)" in the proof of
Theorem~\ref{th:1.1} actually shows
\begin{eqnarray*}
&&{\cal P}(H,\tau)=\bigl\{x(t)=\phi^H_t(x_0)\,\bigm|\, x_0\in L\cap
(\phi^H_{1/2})^{-1}(L)\bigr\}\\
&&\hbox{and so}\quad \sharp{\cal P}(H,\tau)=\sharp
(L\cap\phi^H_{\frac{1}{2}}(L)).
\end{eqnarray*}
 So using the results obtained for the Arnold
conjecture on Lagrangian intersections one may get the estimates of
the lower bound of $\sharp{\cal P}(H,\tau)$ under certain
assumptions. For example, it follows from Theorem~\ref{th:1.1} and
\cite[Theorem H]{FuOOO} that if  $M$ is closed, $L$ is semipositive,
and $L\pitchfork\phi^H_{1/2}(L)$ then
$\sharp{\cal P}(H,\tau)\ge \sum{\rm rank}H_\ast(L;\Z_2)$.}
\end{remark}

Recall that a symplectic manifold $(M, \omega)$ without boundary is
said to be  {\bf geometrically bounded} if there exist a
geometrically bounded Riemannian
  metric $\mu$ on $M$ (i.e., its sectional curvature is bounded above by some constant $K>0$
and injectivity radius $i(M, \mu)>0$) and a $\omega$-compatible
almost complex structure $J$ such that such that
$$
\omega(X, JX)\ge\alpha_0\|X\|^2_\mu\quad {\rm and}\quad |\omega(X,
Y)|\le\beta_0\|X\|_\mu\|Y\|_\mu\quad\forall X, Y\in TM
$$
for some positive constants $\alpha_0$ and $\beta_0$ (cf. \cite{Gr},
\cite{AuLaPo}, \cite{CGK}, \cite{Lu1}). For a real symplectic
manifold $(M,
 \omega,\tau)$ without boundary, if the almost complex structure
 $J$ above can be chosen in $\R{\cal J}(M,\omega)$ we say
 $(M,\omega,\tau)$ to be {\bf real geometrically bounded} (with
 respect to $(J,\mu)$).

The {\bf rationality index} of a symplectic manifold $(M, \omega)$ is defined by
$$
r_\omega=r(M,\omega):=\inf\bigl\{\,\langle \omega, A\rangle\,\bigm|\, A\in
\pi_2(M),\, \langle \omega, A\rangle>0\,\bigr\}\in [0, +\infty],
$$
where we use the convention that the infimum over the
empty set is equal to $+\infty$.
Since $\{\omega(A)\,|\, A\in\pi_2(M)\}$ is a subgroup of $(\mathbb{R}, +)$,
it is easily checked that $r_\omega$ is a
finite positive number if and only if
$\omega\bigl(\pi_2(M)\bigr)=r_\omega\Z$.
 For  $J\in {\cal J}(M, \omega)$ let $m(M, \omega, J)\in [0, +\infty]$
denote  the infimum of the area of all nonconstant $J$-holomorphic
spheres in $M$, where as usual we understand  $m(M, \omega, J)=\infty$ if no nonconstant
$J$-holomorphic sphere exists. Clearly, $r_\omega\le m(M,\omega,J)\;\forall J\in\mathcal{J}(M,\omega)$.
As showed by (ii)-(iii) of Example~\ref{ex:4.2},  there exist
closed symplectic manifolds $(M,\omega)$ such that
$$
0<r_\omega<\sup_{J\in{\cal J}(M,\omega)}m(M,\omega, J)=+\infty.
$$
If $M$ is compact, it directly follows from the Gromov
compactness theorem that
%\begin{equation}\label{e:1.10.1}
$$m(M, \omega, J)>0.$$
%\end{equation}
If $(M, \omega, J)$ is only geometrically bounded as above, this may be derived from
the monotonicity principle (\cite[Prop.4.3.1(ii)]{Sik}): \textsf{For $r_0=\min\{i(M, \mu), \pi/\sqrt{K}\}$,
a compact Riemann surface with boundary $S$ and a $J$-holomorphic map $f:S\to M$,
assume that there exists a $\mu$-metric ball $B(x,r)$ with $r\le r_0$ and with $x\in f(S)$
such that $f(\partial S)\subset\partial B(x,r)$, then}
$$
\frac{\pi\alpha_0}{4\beta_0}r^2\le{\rm Area}_\mu(f(S))\le\frac{1}{\alpha_0}\int_Sf^\ast\omega.
$$
In fact, put $\delta=\min\{\frac{\pi\alpha_0}{4\beta_0}r^2_0,  \frac{\pi\alpha_0}{4\beta_0}(i(M, \mu))^2/9\}$.
It follows that $\int_\Sigma u^\ast\omega\ge \delta$
for every nonconstant $J$-holomorphic map $u$ from a closed Riemann surface $\Sigma$ to $M$.
(See the proof of \cite[Lemma 8.1]{FuO} below Lemma 8.10 therein).

Based on Hofer'method in \cite{Ho2} we can get our second
result.

\begin{theorem}\label{th:1.3}
Let $(M,\omega,\tau)$ be a real geometrical bounded  symplectic
manifold with respect to $J\in\R{\cal J}(M,\omega)$ and a Riemannian
metric $\mu$, and $L={\rm Fix}(\tau)$ be a nonempty compact
submanifold without boundary. Let $H\in C^\infty_0(\R/\mathbb{Z}\times M)$ be a
 symmetric Hamiltonian function. If $r_\omega>0$ and
 $\|H\|<r_\omega$, then
\begin{equation}\label{e:1.11}
\sharp(L\cap{\rm Fix}(\phi^H_1))\ge\sharp {\cal P}_0(H,\tau)\ge {\rm
Cuplength}_{\Z_2}(L),
\end{equation}
and $\sharp {\cal P}_0(H,\tau)\ge {\rm Cuplength}_{\Z}(L)$ if $L$ is
orientable.
\end{theorem}

Note that $r_\omega\in (0, +\infty)$ (resp. $=+\infty$) implies $\omega(\pi_2(M, L))=\frac{r_\omega}{2}\mathbb{Z}$ (resp. $=0$).
As a direct consequence of (\ref{e:1.6}) and Theorems~\ref{th:1.1},
\ref{th:1.3} we get

\begin{theorem}\label{th:1.4}
Let $(M,\omega,\tau, J, L)$ be as in Theorem~\ref{th:1.3}.
If $r_\omega>0$, then every Hamiltonian diffeomorphism $\phi$ on $M$ generated by a
Hamiltonian function $H\in C_0^\infty([0, 1]\times M)$ with
$\|H\|<r_\omega/2$, satisfies the estimates
\begin{equation}\label{e:1.12}
\sharp(L\cap\phi(L))\ge {\rm
Cuplength}_{\Z_2}(L),
\end{equation}
and $\sharp(L\cap \phi(L))\ge {\rm Cuplength}_{\Z}(L)$ if $L$ is
orientable.
\end{theorem}

\begin{remark}\label{th:1.5}
{\rm If $M$ is closed, (\ref{e:1.12}) is a special case of the main result in \cite{Liu}
proved with Floer homology;
the latter and Theorems~\ref{th:1.1}  can only lead to
$\sharp {\cal P}(H,\tau)\ge {\rm Cuplength}_{\Z_2}(L)$,
which is weaker than the second inequality in (\ref{e:1.11}). Actually,
the main result in \cite{Liu} can also be proved by refining
Hofer' arguments in \cite{Ho2} as done in this paper. Hofer'method does not involve
Floer and Morse homologies (and thus complicated transversality arguments).
Recently, Albers and Hein \cite{AH} gave an abstract result based on Morse cohomology.
As in the proof of \cite[Theorem~5.1]{AH}, it may lead to (\ref{e:1.11}),
but no better result.}
\end{remark}

 The twisted product $(\widehat M, \widehat\omega)=(M\times M, \omega\times(-\omega))$
  of a symplectic manifold $(M,\omega)$ and itself with anti-symplectic
involution given by
$$
\tau: M\times M\to M\times M,\;(x,y)\mapsto (y,x),
$$
is a real symplectic manifold with ${\rm Fix}(\tau)=\triangle_M$.
For any $J\in{\cal J}(M,\omega)$ it is easily checked that
$J\times (-J)\in\R{\cal J}(M\times M, \omega\times(-\omega))$ and
\begin{equation}\label{e:1.13}
 m(M\times M, \omega\times(-\omega),
J\times(-J))=2m(M,\omega,J).
\end{equation}
 If  $H\in C^\infty(\R/\mathbb{Z}\times M)$, then
$$
\widehat H:\R\times M\times M\to\R,\;(t, x, y)\mapsto H_t(x)+
H_{-t}(y),
$$
is $1$-periodic in time and symmetric. Note that $X_{\widehat
H_t}(x,y)=(X_{H_t}(x), -X_{H_{-t}}(y))$ by the definition of $X_{H_t}$
above (\ref{e:1.2}). One easily proves that
$z=(x,y):\R/\Z\to\R$ belongs to ${\cal P}(\widehat H,\tau)$ (resp.
${\cal P}_0(\widehat H,\tau)$) if and only if $x\in{\cal P}(H)$
(resp. $x\in{\cal P}_0(H)$) and $y(t)=x(-t)\,\forall t\in\R$.
Here ${\cal P}(H)$ (resp. ${\cal P}_0(H)$) always denote the set of
1-periodic solutions (resp. contractible 1-periodic solutions) of
the equation $\dot x=X_H(t,x)$.
Moreover,
$$
\|\widehat H\|=\int^1_0[\sup_{(x,y)} H_t(x,y)-\inf_{(x,y)}
H_t(x,y)]dt=2\|H\|
$$
and $r_{\widehat{\omega}}=r_\omega$ are clear. Using this and (\ref{e:1.13}) we derive from
Theorem~\ref{th:1.3}:

\begin{theorem}\label{th:1.7}
Let $(M,\omega)$ be a closed  symplectic manifold, and  $H\in C^\infty(\R/\mathbb{Z}\times M)$
satisfy $\|H\|<r_\omega$. Then
 $$
 \sharp{\rm Fix}(\phi^H_1)\ge{\rm
Cuplength}_{\Z_2}(M)\quad\hbox{and}\quad\sharp{\rm Fix}(\phi^H_1)\ge{\rm Cuplength}_{\Z}(M).
$$
\end{theorem}

The first inequality was proved in \cite[Theorem 1.1]{Sch} by Floer homology method.
It is a generalization of the result in \cite{Fl2, Ho2}.
Without the assumption ``$\|H\|<r_\omega$'',  Le and Ono \cite{LeO} got the estimates (\ref{e:1.9}) for $\F=\Z_2$ if $(M,\omega)$ is negative monotone and has
minimal Chern number $N\ge \dim M/2$ (cf. Example~\ref{ex:4.2}(iii) for these two notions).

 The cotangent bundle of a manifold $N$, $(T^\ast N,\omega_{\rm
 can}=-d\lambda_{\rm can})$, is a real symplectic manifold with
the anti-symplectic involution given by
$$
\tau: T^\ast N\to T^\ast N,\;(q,
p)\mapsto (q, -p),
$$
where $q\in N$ and $p\in T^\ast_q N$. Recall that the Liouville
$1$-form $\lambda_{\rm can}$  on $T^\ast N$ is defined by
$\lambda_{\rm can}(\xi)=p(T\pi^\ast\xi)\;\forall\xi\in T_pT^\ast N$,
where $\pi^\ast:T^\ast N\to N$ is the natural projection. The fixed
point set ${\rm Fix}(\tau)$ is the zero section $0_N$ which can be
identified with $N$. Assume now that $N$ is closed. As in \cite{CGK,
Lu1} we can prove that $(T^\ast N,\omega_{\rm can},\tau)$ is
geometrically bounded for some $J\in\R{\cal J}(T^\ast N,\omega_{\rm
can})$ and some metric $G$ on $T^\ast N$.
 Applying Theorem~\ref{th:1.3} to $(T^\ast N,\omega_{\rm can},\tau)$ we
immediately obtain:

\begin{corollary}\label{cor:1.9}
 Let $N$ be a closed manifold, and
$H\in C_0^\infty(\R/\Z\times T^\ast N)$ satisfy  $H(-t, q,p)=H(t, q,
-p)$ for all $t\in\R$ and $(q,p)\in T^\ast N$.
 Then
 $$
\sharp(0_N\cap{\rm Fix}(\phi^H_1))\ge\sharp {\cal P}_0(H,\tau)\ge
{\rm Cuplength}_{\Z_2}(N),
$$
and $\sharp {\cal P}_0(H,\tau)\ge {\rm Cuplength}_{\Z}(N)$ if $N$ is
orientable.
\end{corollary}

%%%%%%%%%%%%%%%%%%%%%%%%%%%%%%%%%%%%%%%%%%%%%%%%%%%%%%%%%%%%%%%%%%%%%%%%
%%Since it was proved in \cite[Theorem 0.4.2]{Cha} that every
%%Hamiltonian diffeomorphism on $(T^\ast N,\omega_{\rm can})$ can be
%%generated by some Hamiltonian $H\in C_0^\infty([0,1]\times T^\ast
%%N)$, Corollary~\ref{cor:1.9}
%%%%%%%%%%%%%%%%%%%%%%%%%%%%%%%%%%%%%%%%%%%%%%%%%%%%%%%%%%%%%%%%%%%%%

This and Theorem~\ref{th:1.1} immediately
lead to

\begin{corollary}\label{cor:1.10}
{\rm ([Ho1, LaSi])}\quad Let $N$ be a closed manifold. Then any
Hamiltonian diffeomorphism $\phi$ on  $(T^\ast N,\omega_{\rm can})$
generated by a Hamiltonian function $H\in C_0^\infty([0, 1]\times T^\ast N)$
satisfies estimates: $\sharp(N\cap\phi(N))\ge {\rm Cuplength}_{\Z_2}(N)$, and
$\sharp(N\cap\phi(N))\ge {\rm Cuplength}_{\Z}(N)$ if $N$ is
orientable.
\end{corollary}

%%%%%%%%%%%%%%%%%%%%%%%%%%%%%%%%%%%%%%%%%%%%%%%%%%%%%%%%%%%%%%%%%%%%%%%%%%%%%%%%%%%%%%%%%%
%%Our following Theorem~\ref{th:1.3} can improve all these results.
%%So far we cannot directly derive the second estimate in Theorem~\ref{th:1.3} from the known
%%results yet.
%%Theorem~\ref{th:1.4} means that the degenerated Arnold-Givental conjecture
%%holds true in the real symplectic manifold $(M, \omega, \tau)$ if $r_\omega=+\infty$.
%%%%%%%%%%%%%%%%%%%%%%%%%%%%%%%%%%%%%%%%%%%%%%%%%%%%%%%%%%%%%%%%%%%%%%%%%%%%%%%%%%%%%%%

 The arrangements of the paper as
follows. In Section~\ref{sec:2.1} we first prove
Theorem~\ref{th:1.1}. Then in Section~\ref{sec:2.2} we complete
the proof of Theorem~\ref{th:1.3} by improving the arguments in
\cite[\S6.4]{HoZe} (also see \cite{Ho2}). Unlike they consider the
space of all bounded trajectories we here only use a subset of it.
 Another different point is
to introduce a definition of topological degree for maps from a
Banach Fredholm bundle to a manifold, not using the $\Z_2$-degree
for Fredholm section having Fredholm index zero as in
\cite[\S6.4]{HoZe}.  The final Section~\ref{sec:3} gives two
examples and a further programme.

\noindent{\bf Acknowledgements}: The results of this paper were
reported in the workshop on Floer Theory and Symplectic Dynamics at
CRM of University of Montreal, May 19-23, 2008.  I would like to
thank the organizers for their invitation, and CRM for hospitality.

\section{Proofs of Theorems~\ref{th:1.1}, \ref{th:1.3}}\label{sec:2}
\setcounter{equation}{0}

\subsection{Proof of Theorem~\ref{th:1.1}}\label{sec:2.1}

\noindent{\bf (i)$\Longrightarrow$ (ii)}:\quad Let $\phi_t$ be
the Hamiltonian flow generated by $H$. Define $Q:[0,1]\times M\to\R$
by $Q(t,x)=H(t/2,x)$, and denote by $\varphi_t$ the flow of $X_Q$.
It is easily proved that
\begin{equation}\label{e:2.1}
\phi_{\frac{1}{2}}=\varphi_1\quad{\rm and}\quad
\|Q\|=\frac{1}{2}\|H\|<\Lambda.
\end{equation}
It follows from (i) that
$$
\sharp(L\cap\phi_{\frac{1}{2}}(L))\ge m.
$$
For any $x_0\in L\cap\phi^{-1}_{\frac{1}{2}}(L)$,
$x(t):=\phi_t(x_0)$ satisfies $\dot x(t)=X_{H_t}(x(t))\;\forall t$
and $x(\frac{1}{2})=\phi_{\frac{1}{2}}(x_0)\in L$. Since
$H_t=H_{1-t}\circ\tau$, for $\frac{1}{2}\le t\le 1$ we have
\begin{eqnarray*}
&&\dot
x(t)=X_{H_t}(x(t))=-d\tau(\tau(x(t)))X_{H_{1-t}}(\tau(x(t)))\quad{\rm
or\; equivalently}\\
&&\frac{d}{dt}\tau(x(t))=-X_{H_{1-t}}(\tau(x(t))).
\end{eqnarray*}
It follows that $y(t):=\tau(x(1-t))$ on $[0, \frac{1}{2}]$
satisfies $\dot y(t)=X_{H_t}(y(t))$. Note that $x(\frac{1}{2})\in L$
implies $y(\frac{1}{2})=\tau(x(\frac{1}{2}))=x(\frac{1}{2})$, i.e.,
$\phi_{\frac{1}{2}}(y(0))=\phi_{\frac{1}{2}}(x_0)$. Hence $y(0)=x_0$ and thus
$\tau(x(1-t))=y(t)=x(t)\;\forall 0\le t\le\frac{1}{2}$.
This implies $x(1-t)=\tau(x(t))\;\forall t\in [0,1]$.
In particular, we get $x(1)=\tau(x_0)=x_0=x(0)$. Moreover, since
$H_0=H_1$, one has $\dot x(1)=\dot x(0)$. Hence $x$ is a
$1$-periodic solution of $\dot x(t)=X_{H_t}(x(t))$ satisfying
$x(1-t)=\tau(x(t))\;\forall t$, that is, $x\in{\cal P}(H,\tau)$. It
is also clear that two different $x_0, x_0^\ast\in
L\cap\phi^{-1}_{\frac{1}{2}}(L)$ give two different elements in ${\cal P}(H,\tau)$,
$x(t)=\phi_t(x_0)$ and $x^\ast(t)=\phi_t^\ast(x_0)$.

 Conversely, each $x\in{\cal P}(H,\tau)$  determines a point $x(0)\in L\cap
\phi^{-1}_{\frac{1}{2}}(L)$ uniquely. So we get
\begin{equation}\label{e:2.2}
{\cal P}(H,\tau)=\{x(t)=\phi_t(x_0)\,|\, x_0\in
L\cap\phi^{-1}_{\frac{1}{2}}(L)\}
\end{equation}
which implies $\sharp{\cal P}(H,\tau)=\sharp
(L\cap\phi_{\frac{1}{2}}(L))\ge m$.\\

\noindent{\bf (ii)$\Longrightarrow$ (i)}:\quad By the assumption
there exists a Hamiltonian $H\in C_0^\infty([0,1]\times M)$ with
$\|H\|<\Lambda$, such that its Hamiltonian flow $\phi_t$
satisfies $\phi_1=\phi$. The proof will be finished along the
line of proof of  \cite[Propsition 2.1.3]{BiPoSa}. Take a small
$\delta>0$ so that $2\|H\|+ 2\delta<2\Lambda$.  Then choose a smooth
function $\lambda:[0, 1]\to [0, 1]$ such that for a given small
$0<\epsilon\ll 1/2$,
\begin{equation}\label{e:2.3}
\left.\begin{array}{lc} \lambda(t)=0\;\hbox{for}\;t\in [0,
\epsilon],\\
 \lambda(t)=0\;\hbox{for}\;t\in [1-\epsilon, 1],\\
\lambda'(t)>0\;\hbox{for}\; t\in (\epsilon, 1-\epsilon).
\end{array}\right\}
\end{equation}
Clearly, $\int^1_0\lambda'(t)dt=1$. Take a time independent
   compactly supported function $F:M\to\R$ which is $\tau$-invariant, such that
   $\|F\|_{C^0}<\delta/4$. Let $f_t$ be the Hamiltonian flow  generated by $F$.
Then the Hamiltonian isotopy
   $\varphi_t := f_{t-\lambda(t)} \circ \phi_{\lambda(t)}$
      is generated by the Hamiltonian function
   $$
   \overline H_t := F + \lambda'(t)(H_{\lambda(t)}-F)\circ f_{\lambda(t)-t}.
   $$
   The function $\overline H_t$ equals $F$ near $t=0$ and $t=1$ and
   hence defines a smooth Hamiltonian on $S^1\times M$.  Moreover,
   $\varphi_1=\phi_1$.
Denote by
$$
A_H(t)=\sup_{x\in M}H_t(x)-\inf_{x\in M}H_t(x)\quad\forall t\in
[0,1].
$$
Then $\|H\|=\int^1_0A_H(t)dt$, and it is easily computed that
\begin{eqnarray*}
A_{\overline H}(t):\!\!\!\!\!&&=\sup_{x\in M}\overline
H_t(x)-\inf_{x\in
M}\overline H_t(x)\\
&&\le \lambda'(t)\bigl(\sup_{x\in M}H_{\lambda(t)}(x)-\inf_{x\in
M}H_{\lambda(t)}(x)\bigr)+ 2\|F\|_{C^0}+ 2\lambda'(t)\|F\|_{C^0}.
\end{eqnarray*}
From this and (\ref{e:2.3}) we arrive at
\begin{eqnarray*}
\|\overline H\|=\int^1_0A_{\overline
H}(t)dt&&\le\int^1_0\lambda'(t)A_H(\lambda(t))dt+ 4\|F\|_{C^0}\\
&&=\int^{1-\epsilon}_\epsilon A_H(\lambda(t))d\lambda(t)+ 4\|F\|_{C^0}\\
&&=\int^1_0 A_H(t)dt+ 4\|F\|_{C^0}\\
&&=\|H\|+ 4\|F\|_{C^0}.
\end{eqnarray*}
Let us define a smooth Hamiltonian $G:[0,1]\times M\to\R$  by
$$
G_t(x)=\left\{
\begin{array}{ll}
2\overline H_{2t}(x)\;&{\rm if}\;0\le t\le 1/2,\\
2\overline H_{2(1-t)}(\tau x)\;&{\rm if}\; 1/2\le t\le 1.
\end{array}\right.$$
It is easy to see that $G_t=F$ near $t=0, 1/2, 1$, and
$G_{1-t}(x)=G_t(\tau x)$ for any $(t,x)\in [0,1]\times M$. Extend
$G$ to $\R\times M$ $1$-periodically in $t$, still denoted by $G$,
we easily see that $G$ satisfies
$\|G\|=2\|\overline H\|<2\|H\|+2\delta<2\Lambda$
and (\ref{e:1.1}), i.e.,
$$
G_{t+1}=G_t\quad{\rm  and}\quad G_{-t}(x)=G_t(\tau x)\;\;\forall
(t,x)\in\R\times M.
$$
It follows that
$$
X_{G_t}(x)=\left\{
\begin{array}{ll}
2X_{\overline H_{2t}}(x)\;&{\rm if}\;0\le t\le 1/2,\\
-2d\tau(\tau x)X_{\overline H_{2(1-t)}}(\tau x)\;&{\rm if}\; 1/2\le
t\le 1
\end{array}\right.
$$
and thus the flow $\psi_t$ of $X_{G}$ and the flow
$\varphi_t$ of $X_{\overline H}$ satisfy
$$
\psi_{t/2}(x)=\varphi_{t}(x)\quad\hbox{for}\;(t,x)\in
[0,1]\times M.
$$
Specially, we have $\psi_{1/2}=\varphi_1=\phi$. Now
for any $y\in{\cal P}(G,\tau)$, the map $x:[0, 1]\to M$ defined by
$x(t)=y(t/2)$ satisfies $\dot x(t)=X_{\overline H_t}(x(t))$. Note
that both $x(0)=y(0)=y(1)$ and $x(1)=y(1/2)$ belong to $L={\rm
Fix}(\tau)$. Hence $x(1)=\varphi_1(x(0))\in
L\cap\varphi_1(L)= L\cap\phi_1(L)=L\cap\phi(L)$
since $\varphi_1=\phi_1=\phi$.

Moreover, for two different $y_1, y_2\in{\cal P}(G,\tau)$
we have $y_1(t_0)\ne y_2(t_0)$ for some $t_0\in [0,1/2]$.
For $t\in [0,1]$ let $x_i(t)=y_i(t/2)$, $i=1,2$. Both satisfy
 $\dot x(t)=X_{\overline H_t}(x(t))$. Since
$x_1(2t_0)\ne x_2(2t_0)$, that is, $\varphi_{2t_0}(x_1(0))\ne
\varphi_{2t_0}(x_2(0))$, we obtain  $x_1(0)\ne x_2(0)$
and thus $x_1(1)\ne x_2(1)$.

In summary, we have proved  $\sharp(L\cap\phi(L))\ge\sharp{\cal P}(G,\tau)$.
Applying Theorem~\ref{th:1.1}(ii) to $G$ we have also that
$\sharp{\cal P}(G,\tau)\ge {\rm Cuplength}_{\Z_2}(L)$.
The desired claim is proved.
 \hfill$\Box$\vspace{2mm}

\subsection{Proof of Theorem~\ref{th:1.3}}\label{sec:2.2}

Let $(M, \omega,\tau)$ be real geometrical bounded for $J\in\R{\cal
J}(M,\omega)$ and a Riemannian metric $\mu$ on $M$.
 By the assumptions of Theorem~\ref{th:1.3}  there exists a compact
subset $K\subset M$  such that
\begin{equation}\label{e:2.4}
{\rm supp}(H_t)\subset K\;\forall t\in\R,\quad L\subset K\quad{\rm
and}\quad\;\bigcup_{x\in{\cal P}_0(H,\tau)}x(\R)\subset K.
\end{equation}
From now on, we assume $(M,g_J)\subset (\R^N,
\langle\cdot,\cdot\rangle)$ by the Nash embedding theorem. Consider
the standard Riemannian sphere $(S^2=\C\cup\{\infty\}, j)$ and the submanifold
of the Banach manifolds $W^{1,p}(S^2, M)$ for a fixed $p>2$,
$$
{\cal B}=\{w\in W^{1,p}(S^2, M)\,|\, w \;\hbox{is contractible}\}.
$$
Let $E_J\rightarrow S^2\times M$ be the vector bundle, whose fiber
over $(z,m)\in S^2\times M$ consists of all linear maps $\phi\colon
T_z S^2\rightarrow T_mM$ such that $J(m)\phi=-\phi\circ j$.  Due to
the inclusion $W^{1,p}(S^2, M)\hookrightarrow C^0(S^2,M)$, for given
$w\in W^{1,p}(S^2, M)$, we can denote by $\bar w\colon
S^2\rightarrow S^2\times M$ the ``graph map'' $\bar w(z)=(z, w(z))$
and write $\bar w^*E_J\rightarrow S^2$ for the pull back bundle.
There exists a natural
  Banach space bundle ${\cal E}\rightarrow \cal B $  whose fiber ${\cal
E}_w=L^{p}(\bar w^*E_J)$ at $w\in {\cal B}$ consists of all $L^p$
sections of the vector bundle
  $\bar w^*E_J\rightarrow  S^2$.
  The nonlinear Cauchy-Riemannian  operator
$\bar \partial_J$,
$$
\bar \partial_J(w)=dw+J\circ dw\circ j,
$$
can be considered as a smooth section of the bundle ${\cal
E}\rightarrow {\cal B}$.

Denote by $Z_T=[-T, T]\times S^1$ for $T>1$. Take a smooth function
$\gamma:\R\to [0,1]$ such that $\gamma(s)=1$ for $s\le -1$,
$\gamma(s)=0$ for $s\ge 0$, and $\gamma'(s)\le 0$ and for $s\in \R$.
Define
\begin{equation}\label{e:2.4.1}
\gamma_T(s)=\left\{\begin{array}{ll}
 1,&\textrm{$s\in [-T+1, T-1]$,}\\
\gamma(s-T),&\textrm{$s\ge T-1$,}\\
 \gamma(-s-T), &\textrm{$s\le -T+1$.}
  \end{array}\right.
 \end{equation}
Then $\gamma'_T(s)\le 0$ for $s\ge T-1$, and $\gamma'_T(s)\ge 0$ for
$s\le-T+1$. Denote by $\nabla$ the Levi-Civita connection with
respect to the metric $\langle\cdot, \cdot\rangle=g_J(\cdot,
\cdot)$. By the definition of $X_{H_t}$
above (\ref{e:1.2}), $\nabla H_t=-JX_{H_t}$.   For  $(z,m)\in
(S^2\setminus\{0,\infty\})\times M$ let us define
\begin{eqnarray*}
h^T_{J}(z,m)\Bigl(\xi\frac{\partial}{\partial x}|_z+
\eta\frac{\partial}{\partial y}|_z\Bigr) \!\!\!\!\!\!\!\!&&=\,
\xi\Bigl(\frac{\gamma_T(s)e^{-2\pi s}\cos(2\pi
t)}{2\pi}\nabla H_t(m)\\
&&\qquad\qquad - \frac{\gamma_T(s)e^{-2\pi
s}\sin(2\pi t)}{2\pi}J(m)\nabla H_{t}(m)\Bigr)\\
&&\,-\, \eta\Bigl(\frac{\gamma_T(s)e^{-2\pi s}\sin(2\pi
t)}{2\pi}\nabla H_t(m)\\
&&\qquad\qquad + \frac{\gamma_T(s)e^{-2\pi s}\cos(2\pi
t)}{2\pi}J(m)\nabla H_{t}(m)\Bigr)
\end{eqnarray*}
 for $\xi,\eta\in\R$ and $z=e^{2\pi(s+it)}\in\C$. It is easily checked that
 $h^T_{J}(z,m)\circ j=-J\circ h^T_{J}(z,m)$, i.e.,
 $h^T_{J}(z,m)\in (E_J)_{(z,m)}$. Note that
\begin{eqnarray*}
&&0<|z|=e^{2\pi s}\le e^{-2\pi (T+1)}\Longleftrightarrow s\in
(-\infty,
-T-1]\Rightarrow \gamma_T(s)=0,\\
&&\infty>|z|=e^{2\pi s}\ge e^{2\pi (T+1)}\Longleftrightarrow s\in
[T+1, +\infty)\Rightarrow \gamma_T(s)=0.
\end{eqnarray*}
Hence we can define $h^T_{J}(0,m)=0, h^T_{J}(\infty,m)=0$ and get a
smooth family of sections $h^T_{J}:S^2\times M\to E_J$, $T>1$. These
give rise to a smooth family of sections of the Banach bundle ${\cal B}\to {\cal E}$,
$g^T_{J}:{\cal B}\to {\cal E}$, $T>1$, where
$$
g^T_{J}(w)(z)=h^T_{J}(z, w(z))\;\;\forall z\in S^2.
$$
For $\lambda\in [0,1]$ we define
\begin{equation}\label{e:2.5}
{\cal F}_{T,\lambda}:{\cal B}\to {\cal E},\;w\mapsto\bar\partial_Jw+
\lambda g^T_J(w).
\end{equation}
Note that $\tau$ and  the standard complex conjugate $c_S$ on $(S^2,
j)$ induce an involution
\begin{equation}\label{e:2.6}
\tau_B:{\cal B}\to {\cal B},\, w\mapsto \tau\circ w\circ c_S^{-1},
\end{equation}
and its lifting involution
\begin{equation}\label{e:2.7}
\tau_E:{\cal E}\to {\cal E},
\end{equation}
 where for $\xi\in{\cal E}_w$, $\tau_E(\xi)\in {\cal E}_{\tau_B(w)}$ is given by
$$
\tau_E\bigl(\xi)(z, \tau_B(w)(z)\bigr)=d\tau(w(\bar
z))\circ\xi(\bar z, w(\bar z))\circ dc_S(z)\;\;\forall z\in S^2.
$$
 Let ${\cal B}^\tau$ be the set of fixed points of $\tau_{\cal B}$. It is
a Banach submanifold in ${\cal B}$, and $w\in{\cal B}$ sits in
${\cal B}^\tau$ if and only if  $w(\bar z)=\tau(w(z))$ for any $z\in
S^2=\C\cup\{\infty\}$. Moreover, the involution $\tau_E$ induces
bundles homomorphisms on ${\cal E}|_{{\cal B}^\tau}$. Denote by
${\cal E}_{+1}$ (resp. ${\cal E}_{-1}$) the eigenspace associated to
the eigenvalue $+1$ (resp. $-1$) of this homomorphism. Then both
${\cal E}_{+1}$ and ${\cal E}_{-1}$ are Banach subbundles of ${\cal
E}|_{{\cal B}^\tau}$, and  ${\cal E}|_{{\cal B}^\tau} = {\cal
E}_{+1} \oplus {\cal E}_{-1}$. Note also that
\begin{equation}\label{e:2.8}
\bar\partial_J\bigl(\tau_B(w)\bigr)=\tau_E\bigl(\bar\partial_J(w)\bigr)\;\;\forall
w\in{\cal B}.
\end{equation}
   So the restriction $\bar\partial_J|_{{\cal B}^\tau}$ gives rise
to a section of the bundle ${\cal E}^+\to {\cal B}^\tau$.

Since $c_S(0)=0$ and $c_S(\infty)=\infty$, we compute
\begin{equation}\label{e:2.9}
g^T_J(\tau_B(w))(z)=h^T_{J}(z, \tau_B(w)(z))=h^T_{J}(z, \tau(w(\bar
z)))\;\;{\rm for}\;\;z\in S^2.
\end{equation}
Note that (\ref{e:1.2}) implies that for $x\in M$,
$$
\nabla H_{-t}(x)=d\tau(\tau(x))\nabla H_t(\tau(x))\quad{\rm
and}\quad d\tau(x)\circ J(x)=-J(\tau(x))\circ d\tau(x).
$$
From the expression of
$h^T_{J}(z,m)\Bigl(\xi\frac{\partial}{\partial x}|_z+
\eta\frac{\partial}{\partial y}|_z\Bigr)$ above one easily checks
$$
h^T_{J}(z, \tau(w(\bar z)))\Bigl(\xi\frac{\partial}{\partial x}|_z+
\eta\frac{\partial}{\partial y}|_z\Bigr)=d\tau(w(\bar z)) h^T_{J}(z,
w(\bar z))\Bigl(\xi\frac{\partial}{\partial x}|_{\bar z}-
\eta\frac{\partial}{\partial y}|_{\bar z}\Bigr),
$$
that is, $h^T_{J}(z, \tau(w(\bar z)))=d\tau(w(\bar z))\circ h^T_{J}(z,
w(\bar z))\circ dc_S(z)$. So (\ref{e:2.7}) and (\ref{e:2.9}) lead to
\begin{equation}\label{e:2.10}
g^T_J(\tau_B(w))=\tau_E(g^T_J(w))\;\;\forall w\in{\cal B}.
\end{equation}
It follows from (\ref{e:2.8}) and (\ref{e:2.10}) that ${\cal
F}_\lambda$ in (\ref{e:2.5}) satisfies
$$
{\cal F}_{T,\lambda}(\tau_B(w))=\tau_E({\cal
F}_{T,\lambda}(w))\;\;\forall w\in{\cal B},
$$
that is, each ${\cal F}_{T,\lambda}$ is equivariant with respect to
the involutions in (\ref{e:2.6}) and (\ref{e:2.7}).
 Hence the restrictions ${\cal F}_{T,\lambda}|_{{\cal
B}^\tau}$ are the sections of the bundle ${\cal E}^+\to{\cal
B}^\tau$. It is easy to prove that all ${\cal F}_{T,\lambda}|_{{\cal
B}^\tau}$ are Fredholm sections of index $n=\dim L$ (by Lemma~\ref{lem:2.4}
the proof of \cite[Prop.6]{Ho2}).
 Define
\begin{eqnarray*}
&&{\cal Z}^\tau_{T,\lambda}:=\{w\in {\cal B}^\tau\,|\, {\cal
F}_{T,\lambda}(w)=0\}\quad{\rm and}\\
&&{\cal Z}^\tau_T:=\{(\lambda, w)\in [0,1]\times{\cal B}^\tau\,|\,
{\cal F}_{T,\lambda}(w)=0\}.
\end{eqnarray*}
The elliptic regularity arguments show that ${\cal
Z}^\tau_{T,\lambda}$ is contained in $C_c^\infty(S^2, M):=\{w\in C^\infty(S^2, M)\,|\, w \;\hbox{is contractible}\}$.
%%%%%%%%%%%%%%%%%%%%%%%%%%%%%%%%%%%%%%%%%%%%%%%%%%%%%%%%%%%%%%%%%%%%%%%%%
%%$\subset C^\infty(S^2, M)$. The same reasoning
%%yields that the zero locus of any smooth perturbation section of
%%${\cal F}_{T,\lambda}$ is contained in $C_c^\infty(S^2, M)$.
%%%%%%%%%%%%%%%%%%%%%%%%%%%%%%%%%%%%%%%%%%%%%%%%%%%%%%%%%%%%%%%%%%%%%%

\begin{lemma}\label{lem:2.1}
For $w\in {\cal Z}^\tau_{T,\lambda}$,  define $u:Z_\infty\to M$  by
$u=w\circ\phi$, where
$$
\phi:Z_\infty=\R\times S^1\to
S^2\setminus\{0,\infty\},\;(s,t)\mapsto e^{2\pi(s+ it)}
$$
 is the biholomorphism. Then $u$
satisfies
\begin{eqnarray}
&&\displaystyle \partial_su(s,t)+ J(u(s,t))\bigl(\partial_t u(s,t)
-\lambda\gamma_T(s) X_{H_t}(u(s,t)\bigr)=0,\label{e:2.11}\\
&& E(u):=\int_{Z_\infty}|\partial_s u|_{g_J}^2ds dt \le \|H\|\le
2\|H\|_{C^0}.\label{e:2.12}
\end{eqnarray}
\end{lemma}

\noindent{\bf Proof.}
The equation $dw(z)+ J(w)\circ dw(z)\circ j+
h^\lambda_{J,T}(z,w(z))=0$ yields
\begin{eqnarray*}
&&dw(z)(\frac{\partial}{\partial x})+ J(w)\circ dw(z)\circ
j(\frac{\partial}{\partial x})+
h^\lambda_{J,T}(z,w(z))(\frac{\partial}{\partial x})=0,\qquad{\rm that\; is}\\
&&\partial_xw+ J(w)\partial_yw + \frac{\lambda\gamma_T(s)e^{-2\pi
s}\cos(2\pi
t)}{2\pi}\nabla^JH_t(w)\\
&& - \frac{\lambda\gamma_T(s)e^{-2\pi s}\sin(2\pi
t)}{2\pi}J(w)\nabla^JH_{t}(w)=0.
\end{eqnarray*}
Since
$$
\partial_xw+ J(w)\partial_yw=\frac{e^{-2\pi s}\cos(2\pi
t)}{2\pi}(\partial_su+ J(u)\partial_tu)- \frac{e^{-2\pi s}\sin(2\pi
t)}{2\pi}J(u)(\partial_su+ J(u)\partial_tu),
$$
it follows that $u(s,t)=w(e^{2\pi(s+it)})$ satisfies
\begin{eqnarray*}
&&\frac{e^{-2\pi s}\cos(2\pi t)}{2\pi}(\partial_su+
J(u)\partial_tu)- \frac{e^{-2\pi s}\sin(2\pi
t)}{2\pi}J(u)(\partial_su+
J(u)\partial_tu)\\
&&+\frac{\lambda\gamma_T(s)e^{-2\pi s}\cos(2\pi
t)}{2\pi}\nabla^JH_t(u) - \frac{\lambda\gamma_T(s)e^{-2\pi
s}\sin(2\pi t)}{2\pi}J(u)\nabla^JH_{t}(u)\\
&&=\frac{e^{-2\pi s}\cos(2\pi t)}{2\pi}\bigl(\partial_su+
J(u)\partial_tu +
\lambda\gamma_T(s)\nabla^JH_{t}(u)\bigr)\\
&&+\frac{e^{-2\pi s}\sin(2\pi t)}{2\pi}J(u)\bigl(\partial_su+
J(u)\partial_tu + \lambda\gamma_T(s)\nabla^JH_{t}(u)\bigr)=0.
\end{eqnarray*}
This is equivalent to (\ref{e:2.11}) since $\nabla H_t=-JX_{H_t}$ and $g_J(X, JX)=0$ for any
$X\in TM$.

As to (\ref{e:2.12}), note that the contractility of $w:S^2\to M$
implies
\begin{eqnarray*}
&&0=\int_{S^2}w^\ast\omega=\int_{Z_\infty}u^\ast\omega=\int_{Z_\infty}
\bigl(|\partial_su|^2_{g_J}+\lambda\gamma_T(s)dH_t(\partial_su)\bigr)dsdt\\
&&\quad =\int_{Z_\infty} |\partial_su|^2_{g_J}ds
dt+\lambda\int^1_0dt\int^{T+1}_{-T-1}\gamma_T(s)\frac{d}{ds}H_t(u))ds\\
&&\quad =\int_{Z_\infty} |\partial_su|^2_{g_J}ds
dt-\lambda\int^1_0dt\int^{T+1}_{-T-1}\gamma'_T(s)H_t(u))ds.
\end{eqnarray*}
Hence
\begin{eqnarray*}
E(u)\!\!\!\!\!\!\!\!&&=\int_{Z_\infty} |\partial_su|^2_{g_J}ds dt\\
&&=\lambda\int^1_0dt\int^T_{T-1}\gamma'_T(s)H_t(u(s))ds+
\lambda\int^1_0dt\int^{-T+1}_{-T}\gamma'_T(s)H_t(u(s))ds\\
&&\le\lambda\int^1_0\sup_pH_t(p)dt\int^{-T+1}_{-T}\gamma'_T(s)ds
+\lambda\int^1_0\inf_pH_t(p)dt\int^T_{T-1}\gamma'_T(s)ds\\
&&=\lambda\int^1_0\sup_pH_t(p)dt
-\lambda\int^1_0\inf_pH_t(p)dt\le\lambda\|H\|\le 2\|H\|_{C^0},
\end{eqnarray*}
where  the first inequality is because $\gamma'_T(s)\ge 0$ for $-T\le s\le -T+1$,
and $\gamma'_T(s)\le 0$ for $T-1\le s\le T$.
\hfill$\Box$\vspace{2mm}

\begin{lemma}\label{lem:2.2}
Suppose that $\|H\|<+\infty$. Then there exists a compact subset
$W\subset M$ such that $w(S^2)\subset W$ for any $(\lambda, w)\in
{\cal Z}^\tau_T$, and this $W$ can be assumed to be a compact
submanifold of codimension zero and to contain $K$ in its interior.
\end{lemma}

\noindent{\bf Proof.}\quad Define $\Delta(w):=w^{-1}(M\setminus
K)\subset S^2$. As in Lemma~\ref{lem:2.1}, let $u:Z_\infty\to M$ be
defined by $u=w\circ\phi$. By (\ref{e:2.12}) we may derive
$$
\int_{\Delta(w)}w^\ast\omega\le E(u)\le \|H\|.
$$
Then one can complete the proof as in the proof of \cite[Theorem
2.9]{Lu1} or as in the proof of Lemma~\ref{lem:2.3}(i) below.
There exists also another method to prove this.
Each $(\lambda, w)\in {\cal Z}^\tau_T$ satisfies
$\bar\partial_Jw(z)+\lambda h^T_{J}(z, w(z))=0\;{\rm for}\;z\in S^2$.
Thus the ``graph map'' $\bar w\colon S^2\rightarrow S^2\times M$ given by $\bar w(z)=(z, w(z))$
is holomorphic with respect to the almost complex structure $J_{H}$ on
$S^2\times M$  by
$$
J_{H,\lambda}(z, m)(X_1,X_2)=(iX_1, -\lambda J(m)\circ h^T_{J}(z, m)X_1+  J(m)X_2).
$$
Then fixing a metric $\tau$ on $S^2$ and applying \cite[Prop.4.4.1]{Sik} to $\bar w\colon S^2\rightarrow (S^2\times M, J_{H,\lambda}, \tau_0\oplus\mu)$,
the desired conclusion can be obtained.
\hfill$\Box$\vspace{2mm}

Let $C_c^\infty(S^1, M)$ denote the set of all contractible smooth
loops $x:S^1\to M$, and
$$
{\cal L}(M,\tau):=\{x\in C_c^\infty(S^1, M)\,|\,
x(-t)=\tau(x(t))\;\forall t\in\R\}.
$$
In the following we always assume that
$C^\infty(\R\times S^1, M)$ is equipped with the compact open
$C^\infty$-topology. Then it is not necessarily path connected even if
$M$ is so.
For $u\in C^\infty(\R\times S^1, M)$ and
$s\in\R$ we write $u(s):S^1\to M$ by $u(s)(t):=u(s,t)$. It is clear
that $u(s)\in C_c^\infty(S^1, M)\;\forall s\in\R$ if and only if
$u(s)\in C_c^\infty(S^1, M)$ for some $s\in\R$. When
$u\in C^\infty(\R\times S^1, M)$ satisfies the equation
\begin{equation}\label{e:2.13}
 \partial_su(s,t)+ J(u(s,t))(\partial_t u(s,t)
-X_{H_t}(u(s,t))=0,
\end{equation}
we define its energy by $E(u)=\int_{Z_\infty}|\partial_s u|^2_{g_J}ds
dt<+\infty$.
 Denote by
\begin{eqnarray}
&{\cal C}^\tau:=\{u\in C^\infty(\R\times S^1, M)\,|\, u(s)\in{\cal
L}(M,\tau)\;\forall s\in\R\},\label{e:2.14}\\
&X_\infty^\tau:=\Bigl\{u\in{\cal C}^\tau\,\Bigm|\, u\;{\rm
satisfies}\;(\ref{e:2.13}),\;E(u)\le \|H\|\Bigr\}.\label{e:2.15}
\end{eqnarray}
Both are equipped with the topology induced from $C^\infty(\R\times
S^1, M)$.

\begin{lemma}\label{lem:2.3}
{\rm (i)}  The compact submanifold $W$ in
Lemma~\ref{lem:2.2} can be enlarged so that $u(\R\times S^1)\subset
W$ for all $u\in X_\infty^\tau$. \\
{\rm (ii)}  $X^\tau_\infty$ is a compact metrisable space provided that
$\|H\|<m(M,\omega,J)$.\\
{\rm (iii)}
If $\sharp {\cal P}_0(H,\tau)$ is finite, then
for every  $u\in{\cal C}^\tau$
satisfying (\ref{e:2.13}) and $E(u)<+\infty$ there exist $x^+, x^-\in{\cal P}_0(H,\tau)$
such that
$$
\lim_{s\to\pm\infty}u(s,t)=x^\pm(t)\quad\hbox{and}\quad \lim_{s\to\pm\infty}\partial_su(s,t)=0,
$$
where both limits are uniform in the $t$-variable.
\end{lemma}

\noindent{\bf Proof.}\quad (i) We may assume that $M$ is noncompact. Let $u\in{\cal C}^\tau$
satisfy (\ref{e:2.13}) and $E(u)<+\infty$. Then
$$
\int^{+\infty}_{-\infty}\left(\int_{S^1}|\partial_tu(s,t)-X_{H_t}(u(s,t))|^2_{g_J}dt\right)ds=
E(u)<+\infty.
$$
Hence there exist sequences $s_k^+\uparrow +\infty$ and $s_k^-\downarrow-\infty$ such that
\begin{equation}\label{e:2.15.1}
\lim_{k\to+\infty}\left\|\frac{\partial u}{\partial t}(s_k^\pm,\cdot)-X_{H_t}(u(s_k^\pm,\cdot))\right\|^2_{L^2}=0.
\end{equation}
Clearly, we may assume $0<s_1^+<s_2^+<\cdots$ and $0>s_1^->s_2^->\cdots$.
Since $X_{H_t}$ vanishes outside the compact subset $K$, it follows from (\ref{e:2.15.1}) that
there exists a constant $C>0$ such that
\begin{equation}\label{e:2.15.2}
\left\|\frac{\partial u}{\partial t}(s_k^\pm,\cdot)\right\|^2_{L^2}\le C,\quad\forall k=1,2,\cdots.
\end{equation}
These imply that for all $t\in [0,1]$,
\begin{eqnarray*}
d_{g_J}(u(s_k^\pm, t), u(s_k^\pm,0))&\le&\int^t_0\left|\frac{\partial u}{\partial t}(s_k^\pm,\tau)\right|_{g_J} d\tau\\
&\le&\sqrt{t}\left(\int^t_0\left|\frac{\partial u}{\partial t}(s_k^\pm,\tau)\right|^2_{g_J} d\tau\right)^{1/2}\\
&\le & \sqrt{C},\quad\forall k=1,2,\cdots.
\end{eqnarray*}
Since $u(s_k^\pm,0)\in L$, it follows that all $u_k(\{s_k^\pm\}\times S^1)$ are contained
in a compact subset $\hat{K}$ of $M$. Clearly, we can assume that
$\hat{K}$ is a compact submanifold of codimension zero and with boundary and that
$K$ is contained the interior of $\hat{K}$.

Now let us assume that this $u$ belongs to $X_\infty^\tau$. Define
$$
w=u|_{[s^-_1, s^+_1]\times S^1},\quad w_k^+=u|_{[s^+_k, s^+_{k+1}]\times S^1},\quad
w_k^-=u|_{[s^-_{k+1}, s^-_{k}]\times S^1},\quad k=1,2,\cdots.
$$
Then each connected component $\Sigma$ of
$$
w^{-1}(M\setminus \hat{K})\;\;\hbox{or}\;\;(w_j^+)^{-1}(M\setminus \hat{K})\;\;\hbox{or}\;\;
(w_j^-)^{-1}(M\setminus \hat{K})\;\hbox{($j\in\mathbb{N}$)},
 $$
 has compact closure and  is the increasing union of connected compact Riemannian surfaces
with smooth boundary  $\Sigma_j$, $j=1,2,\cdots$.
For sufficiently large $j$ we have always $\partial\Sigma_j\subset\hat{K}_1=\{x\in M\,|\,
d(x,\hat{K})\le 1\}$. Note that
the restriction of $w$ or $w_j^+$ or $w_j^-$ to each $\Sigma_j$ is $J$-holomorphic
and has the energy $\le \|H\|$. By \cite[Prop.4.4.1]{Sik} we may deduce that
this restriction has the image contained in the $\tau$-neighborhood of $(\hat{K_1})_\delta$
for some $\tau>0$ only depending on $(M, \omega, \mu, J)$.\\

\noindent{(ii)} By (i) we may assume that $M$ is compact below.
As in \cite[page 236]{HoZe}, it
suffices to prove that there exists a constant $C>0$ such that
\begin{equation}\label{e:2.16}
 |\nabla u(s,t)|_{g_J}\le C\quad\forall
u\in X^\tau_\infty\;{\rm and}\;(s,t)\in Z_{\infty}.
\end{equation}
Arguing indirectly, as on pages 236-238 in \cite{HoZe}, we find
sequences $\varepsilon_k\downarrow 0$, $\{t_k\}_k\subset [0,1]$ and
$\{u_k\}_k\subset X^\tau_\infty$ such that
$$
\left.\begin{array}{ll} t_k\to t_0\in
 [0,1],\; \varepsilon_kR_k\to +\infty\quad{\rm for}\; R_k=|\nabla u_k(0,t_k)|_{g_J}\to +\infty,\\
 |\nabla u_k(s,t)|_{g_J}\le 2 |\nabla u_k(0,t_k)|_{g_J}\quad{\rm
 if}\;|s|^2+ |t-t_k|^2\le\varepsilon_k^2,\; 0\le t_k\le 1
 \end{array}\right\}
 $$
where we consider the $u_k$ as maps defined on $\R\times\R$ by a
1-periodic continuation in the $t$-variable. It follows that the new
sequence $v_k\in C^\infty(\R^2, M)$ defined by
$$
v_k(s,t)=u_k\biggl(\frac{s}{R_k}, t_k+\frac{t}{R_k}\biggr)\;{\rm
for}\; s^2+ t^2\le (\varepsilon_kR_k)^2
$$
converges, in $C^\infty(\R^2, M)$, to $v\in C^\infty(\R^2, M)$ which
satisfies
\begin{equation}\label{e:2.17}
|\nabla
v(0)|_{g_J}=1,\;\sup_{x\in\R^2}|\nabla v(x)|_{g_J}\le 2,\;v_s+
J(v)v_t=0.
\end{equation}
Denote by $B(p, r)\subset\R^2$ the disk centred at $p$ and of radius
$r$. Then
\begin{eqnarray*}
\int_{B(0,\varepsilon_kR_k)}|\partial_sv_k|_{g_J}^2dsdt &&=
\int_{B(0,\varepsilon_kR_k)}\frac{1}{R_k^2}\biggl|\partial_su_k(\frac{s}{R_k},
t_k+ \frac{t}{R_k})\biggr|_{g_J}^2dsdt\\
&&=\int_{B((0, t_k), \varepsilon_k)}\bigl|\partial_su_k(s,
t)\bigr|_{g_J}^2dsdt\\
&&\le E(u_k)\le\|H\|
\end{eqnarray*}
for sufficiently large $k$ (so that $\varepsilon_k<1/2$). It easily
follows that
$$
\int_{\C}|\partial_s v|^2_{g_J}ds dt\le\|H\|<m(M,\omega,J).
$$
However, (\ref{e:2.17}) and Gromov's removable singularity allow us
to extend $v$ to a nonconstant $J$-holomorphic sphere $v_\infty:
S^2\to M$ with
$$
\int_{S^2}v^\ast_\infty=\int_{\C}|\partial_s v|^2_{g_J}ds
dt\le\|H\|<m(M,\omega,J)
$$
which contradicts to the definition of $m(M,\omega,J)$.
(\ref{e:2.16}) is proved.

\noindent{(iii)}. Since the condition $\sharp\mathcal{P}(H)<+\infty$
is actually sufficient for ``(i)$\Rightarrow$(ii)" in the proof of \cite[Prop.1.21]{Sa},
we may complete the proof with the same reason.
 \hfill$\Box$\vspace{2mm}

\begin{lemma}\label{lem:2.4}
Suppose that  $\|H\|<m(M,\omega,J)$. Then ${\cal Z}^\tau_{T,\lambda}$
 and ${\cal Z}^\tau_T$ are compact
in $C^\infty(S^2, M)$ and $[0,1]\times C^\infty(S^2, M)$, respectively.
\end{lemma}

\noindent{\bf Proof}. By Lemma~\ref{lem:2.2} we may assume $M$ to be
compact. Using (\ref{e:2.12}) we can, as in the proof of
Lemma~\ref{lem:2.3}, prove that there exists a constant $C_T>0$ such
that for every $(\lambda, w)\in{\cal Z}^\tau_T$ and
$u=w\circ\phi:Z_\infty\to M$ as in Lemma~\ref{lem:2.1},
\begin{equation}\label{e:2.18}
 \sup_{(s,t)\in Z_\infty}|\nabla u(s,t)|_{g_J}\le C_T.
\end{equation}
It implies that for each multi-index $\alpha\in\N^2$ one can find a
constant $C_{T,\alpha}>0$ such that for all $u$ as above,
\begin{equation}\label{e:2.19}
 \sup_{(s,t)\in Z_\infty}|(D^\alpha u)(s,t)|_{g_J}\le C_{T,\alpha}.
\end{equation}

Now suppose that ${\cal Z}^\tau_T$ is noncompact. Then there exists
 sequences $\{(\lambda_k, w_k)\}_k\subset{\cal Z}^\tau_T$ and
 $\{z_k\}_k\subset S^2=\CP^1$ such that
\begin{equation}\label{e:2.19.1}
 \lambda_k\to\lambda_0\quad{\rm and}\quad |dw_k(z_k)|=\|dw_k\|:=\max_{z\in S^2}|dw_k(z)|\to+\infty,
\end{equation}
where $|dw_k(z)|$ is the norm of the tangent map $dw_k(z):T_zS^2\to
T_{w_k(z)}M$ induced by $g_J$ and the standard Riemannian metric on
$S^2$. We may assume that $z_k\to z_0\in S^2=\CP^1$. By
(\ref{e:2.18}) this $z_0$ must be $0$ or $\infty$ in $\CP^1$.
(Otherwise, passing to a subsequence we may assume
$\inf_k d(z_k, 0)>0$ and $\inf_k d(z_k, \infty)>0$.
Thus there exists a large $T_0>0$ such that $(s_k,t_k)=\phi^{-1}(z_k)\subset S^1\times [-T_0, T_0]$
for all $k$. It follows from (\ref{e:2.19.1}) that $u_k=w_k\circ\phi$
satisfies $|du_k(s_k,t_k)|\to\infty$, which contradicts to (\ref{e:2.19}).)
By the
Gromov compactness theorem the sequence $\{w_k\}_k$ has a subsequence,
still denoted by $\{w_k\}_k$, converges weakly to a connected union of
$N\ge 1$ nonconstant $J$-holomorphic spheres $v_1,\cdots, v_N:S^2\to
M$ and a smooth map $w_\infty:S^2=\CP^1\to M$ satisfying
\begin{equation}\label{e:2.20}
\bar\partial_Jw+ \lambda_0 g^T_J(w)=0.
\end{equation}
In particular, $[v_1\sharp\cdots\sharp v_N\sharp
w_\infty]=0\in\pi_2(M)$. Let $u_\infty=w_\infty\circ\phi:Z_\infty\to
M$. Then as in the proof of Lemma~\ref{lem:2.1} we have
\begin{eqnarray*}
&&0=\sum^N_{k=1}\int_{S^2}v_k^\ast\omega+
\int_{S^2}w_\infty^\ast\omega=\sum^N_{k=1}\int_{S^2}v_k^\ast\omega+\int_{Z_\infty}u_\infty^\ast\omega\\
&&\; =\sum^N_{k=1}\int_{S^2}v_k^\ast\omega+\int_{Z_\infty}
\bigl(|\partial_su_\infty|^2_{g_J}+\lambda_0\gamma_T(s)dH_t(\partial_su_\infty)\bigr)dsdt\\
&&\;=\sum^N_{k=1}\int_{S^2}v_k^\ast\omega+\int_{Z_\infty}
|\partial_su_\infty|^2_{g_J}ds
dt+\lambda_0\int^1_0dt\int^{T+1}_{-T-1}\gamma_T(s)\frac{d}{ds}H_t(u_\infty))ds.
\end{eqnarray*}
It follows that
\begin{eqnarray*}
&&m(M,\omega,J)\le Nm(M,\omega, J)+
E(u_\infty)\\
&&\le \sum^N_{k=1}\int_{S^2}v_k^\ast\omega + E(u_\infty)\\
&&=-\lambda_0\int^1_0dt\int^{T+1}_{-T-1}\gamma_T(s)\frac{d}{ds}H_t(u_\infty))ds\\
&&\le\lambda_0\|H\|\le\|H\|<m(M,\omega,J)
\end{eqnarray*}
as in the proof of Lemma~\ref{lem:2.1}.
This contradiction gives the desired conclusion.
\hfill$\Box$\vspace{2mm}

For $T>1$  we set
\begin{eqnarray*}
&&X^\tau_T:=\{u\in C^\infty(Z_T, M)\,|\, u(0)\in{\cal
L}(M,\tau)\;{\rm and}\; \int_{Z_T}|\partial_su|^2_{g_J}\le \|H\| \},\\
&&X^{\tau, J}_T:=\{u\in X^\tau_{T}\,|\,
\partial_su+ J(u)\partial_tu+ \nabla H_t(u)=0\,\;{\rm
on}\;Z_T\}.
\end{eqnarray*}
As in the proofs of (i)-(ii) of Lemma~\ref{lem:2.3} we may get

\begin{lemma}\label{lem:2.5}
The compact submanifold $W$ in Lemma~\ref{lem:2.3} can be
furthermore enlarged so that $u(Z_T)\subset W$ for all $u\in
X_{T}^\tau$. Moreover, there exists a constant $\widetilde C>0$ such
that for every $T>2$,
\begin{equation}\label{e:2.21}
\sup\left\{|\nabla u(s,t)|_{g_J}\,\Big|\, (s,t)\in Z_{T-2}\right\}\le \widetilde C\quad\forall
u\in X^\tau_T.
\end{equation}
\end{lemma}

Let $\gamma_T(s)$ be as in (\ref{e:2.4.1}). Define
\begin{equation}\label{e:2.22}
\sigma_T:X^\tau_{T}\to{\cal C}^\tau,\;u\mapsto\sigma_T(u)
\end{equation}
by $\sigma_T(u)(s,t)=u(\gamma_T(s)s,t)$. Then
$\sigma_T(u)(s,t)=u(s,t)\;\forall (s,t)\in Z_{T+1}$.

\begin{theorem}\label{th:2.6}
Suppose that  $\|H\|<m(M,\omega,J)$. Then for a
given open neighborhood $U$ of $X^\tau_\infty$ in ${\cal C}^\tau$
there exists $T_0>1$ such that
$$
\sigma_T(X^{\tau, J}_{T})\subset U\;\hbox{for any}\;T\ge T_0.
$$
Furthermore, this $T_0$  can be enlarged so that
$$
\sigma_T(u|_{Z_T})\in U\quad\forall T>T_0
$$
for any  $u=w\circ\phi$ with $w\in{\cal Z}^\tau_{T,1}$, where ${\cal Z}^\tau_{T,1}$ is as above
Lemma~\ref{lem:2.1}.
\end{theorem}

\noindent{\bf Proof}. Since (\ref{e:2.21}) implies that for each
multi-index $\alpha\in\N^2$ one can find a constant $\widetilde
C_{\alpha}>0$ such that for every $T>6$,
\begin{equation}\label{e:2.23}
 \sup\left\{|(D^\alpha u)(s,t)|_{g_J}\,\Big|\, (s,t)\in Z_{T-3}\right\}\le \widetilde C_{\alpha}\quad\forall u\in X^\tau_T.
\end{equation}
As in the arguments on pages 244-245 of \cite{HoZe}, suppose that
there exist an open neighborhood $U$ of $X^\tau_\infty$ in ${\cal
C}^\tau$ and sequences $T_k\to+\infty$ and $u_k\in X^\tau_{T_k}$
such that $u_k\notin U$ for all $k$. From (\ref{e:2.23}) we may
choose a subsequence $\{u_{k_j}\}_j$ of $\{u_k\}_k$ such that $u_{k_j}$
converges to $u$ in $C^\infty_{loc}(\R\times S^1, M)$. Clearly, $u$ satisfies
\begin{eqnarray*}
&&\partial_su+ J(u)\partial_tu+ \nabla H_t(u)=0\,\;{\rm
on}\;Z_\infty,\\
&&u(0,\cdot)\in C_c^\infty(S^1, M)\quad{\rm and}\quad u(s,-t)=\tau(u(s,t))\;\forall (s,t)\in Z_\infty,\\
&&E(u)=\int_{Z_\infty}|\partial_su(s,t)|^2_{g_J}dsdt\le\|H\|.
\end{eqnarray*}
That is, $u\in X^\tau_\infty$. Moreover, all $u_{k_j}$ belong to the closed subset ${\cal
C}^\tau\setminus U$ of ${\cal C}^\tau$. Hence $u\notin U$, which contradicts
$u\in X^\tau_\infty\subset U$.
\hfill$\Box$\vspace{2mm}

For ${\cal C}^\tau$ in (\ref{e:2.14}) we define an evaluation map
\begin{equation}\label{e:2.24}
\pi:{\cal C}^\tau\to L,\;u\mapsto u(0,0),
\end{equation}
and denote $\check{H}^\ast$ by the Alexander-Spanier cohomology.
Then Theorem~\ref{th:1.3} can be derived from the following result.

\begin{theorem}\label{th:2.7}
Under the assumptions, for every open neighborhood $U$ of
$X^\tau_\infty$ in ${\cal C}^\tau$ the restriction $\pi|_U$ induces
an injection
$$
(\pi|_U)^\ast: \check{H}^\ast(L,\Z_2)\to \check{H}^\ast(U,\Z_2).
$$
So the continuity property of the Alexander-Spanier cohomology
implies
$$
\pi|_{X^\tau_\infty}:  \check{H}^\ast(L,\Z_2)\to
\check{H}^\ast(X^\tau_\infty,\Z_2)
$$
is injective. If $L$ is orientable, $\pi|_{X^\tau_\infty}:
\check{H}^\ast(L,\Z)\to \check{H}^\ast(X^\tau_\infty,\Z)$ is also
injective.
\end{theorem}

Consequently, ${\rm Cuplength}_{\Z_2}(X^\tau_\infty)\ge {\rm Cuplength}_{\Z_2}(L)$,
 and ${\rm Cuplength}_{\Z}(X^\tau_\infty)\ge {\rm Cuplength}_{\Z}(L)$
 if $L$ is orientable.\\

 \noindent{\bf Proof of Theorem~\ref{th:1.3}}.\quad
Clearly, we may assume ${\cal P}_0(H,\tau)$ to be a finite set under the assumptions
of Theorem~\ref{th:1.3}. Consider the closed $1$-form $\alpha$ on ${\cal L}(M)$ given by
\begin{equation}\label{e:2.25}
{\alpha}_x(\xi)=\int^1_0\omega(\dot x(t)-X_{H_t}(x(t),\xi(t))dt\quad\forall (x,\xi)\in T{\cal L}(M,\tau).
\end{equation}
It restricts to a closed $1$-form $\alpha^\tau$ on ${\cal L}(M,\tau)$.
Let $\phi_\omega:\pi_2(M)\to\mathbb{R}$ be the homomorphism defined by integration of $\omega$.
Denote by $\tilde{\cal L}(M)$ the set of all pairs $\tilde{x}=(x, [w])$, where $x\in {\cal L}(M)$
and $[w]$ is an equivalence class of smooth discs $w:D^2\to M$ with $w(e^{2\pi it})=x(t)\;\forall t$
for the equivalence relation $\sim$: $w\sim w'$ if and only if the sphere $w\sharp\bar{w'}$ being vanished
by $\omega$. Then $\Pi:\tilde{\cal L}(M)\to {\cal L}(M),\;(x,[w])\mapsto x$
is a covering whose desk group is the quotient $\Gamma(\omega)=\pi_2(M)/{\rm ker}(\phi_\omega)$.
The symplectic action functional
\begin{equation}\label{e:2.26}
{\cal A}_{H}:\tilde{\cal L}(M)\to \R,\;(x,[w])\mapsto -\int_{D^{2}}w^{*}\omega +\int_{0}^{1}H(t,x(t))\,dt
\end{equation}
is a primitive of $\Pi^\ast\alpha$, i.e., $d{\cal A}_{H}(x,[w])[\xi]=(\Pi^\ast{\alpha})_{(x,[w])}(\xi)=
\alpha_x(\xi)$
for any $\xi\in C^\infty(x^\ast TM)$. Let ${\cal A}_{H,\tau}$ be the restriction of
${\cal A}_{H}$ to $\tilde{\cal L}(M,\tau):=\Pi^{-1}({\cal L}(M,\tau))$. Then
$d{\cal A}_{H,\tau}=\alpha^\tau$.

By the assumption the rationality index $r_\omega$ of $(M,\omega)$
is positive. If $r_\omega=+\infty$, i.e., $\omega|_{\pi_2(M)}=0$,
then $\tilde{\cal L}(M)={\cal L}(M)$. If $r_\omega\in (0, +\infty)$,
${\cal A}_{H,\tau}$ descends to a map ${\cal A}_{H,\tau}^\ast:{\cal L}(M,\tau)\to\mathbb{R}/r_\omega\mathbb{Z}$,
which is a primitive of $\alpha^\tau$. Let $p:\mathbb{R}\to\mathbb{R}/r_\omega\mathbb{Z}$
be the canonical projection in the latter case. Define $a_H: {\cal L}(M,\tau)\to \R$ by
\begin{equation}\label{e:2.26.1}
a_H(x)=\left\{\begin{array}{ll}
&{\cal A}_{H,\tau}(x)\quad\hbox{if}\:r_\omega=+\infty,\\
 &(p|_{[0, r_\omega)})^{-1}\circ {\cal A}_{H,\tau}^\ast(x)\quad\hbox{if}\:r_\omega\in (0,+\infty).
 \end{array}\right.
\end{equation}
This is continuous and satisfies
\begin{equation}\label{e:2.26.2}
\frac{d}{ds}a_H(u(s))=\int^1_0|\partial_su(s,t)|^2_{g_J}dt\quad\forall u\in X^\tau_\infty,
\end{equation}
where $u(s)(t)=u(s,t)$. By Lemma 2 on \cite[page 225]{HoZe} (or its proof)
$$
\frac{d}{ds}a_H(u(s))|_{s=s_0}=0\quad\hbox{for some $s_0\in\mathbb{R}$}
$$
implies that $u(s)=u(s_0)\;\forall s\in\mathbb{R}$ and
$x:=u(s_0)=u(s_0,\cdot)$ belongs to ${\cal P}_0(H,\tau)$.
This shows that  the natural flow on the compact metric space $X^\tau_\infty$ defined by
\begin{equation}\label{e:2.26.3}
\Phi:\R\times X^\tau_\infty\to X^\tau_\infty:(\sigma, u)\mapsto \sigma\cdot
u,
\end{equation}
where $(\sigma\cdot u)(s,t)=u(\sigma+ s,t)$, is  gradient-like and
has $a_H$ as a Ljapunov function.  Thus Corollary on \cite[page 42]{CoZe}
yields
\begin{equation}\label{e:2.26.4}
\sharp{\cal P}_0(H,\tau)\ge {\rm
Cuplength}_{\Z_2}(L)\quad\hbox{(or $\ge {\rm Cuplength}_{\Z}(L)$ if $L$ is
orientable).}
\end{equation}
This and Theorem~\ref{th:2.7} give the desired conclusion immediately.
\hfill$\Box$\vspace{2mm}

\section{Proof of Theorem~\ref{th:2.7}}\label{sec:3}
\setcounter{equation}{0}

In order to prove this result let us recall that  a {\bf Banach
Fredholm bundle} of index $r$ and with compact zero sets is a triple
$(X, E, S)$ consisting of  a Banach manifold $X$,  a Banach vector
bundle $E\to X$ and  a Fredholm section $S$ of index $r$ and with
compact zero sets. If the determinant bundle ${\rm det}(S)\to Z(S)$
is oriented, i.e., it is trivializable and is given a continuous
section nowhere zero, we said $(X, E, S)$  to be {\bf oriented}. One
has the following standard result (cf. \cite[Theorem 1.5]{LuT}).

\begin{theorem}\label{th:2.9}
Let $(X, E, S)$ be a Banach Fredholm bundle of index $r$.  Then
there exist finitely many smooth sections $\sigma_1,
\sigma_2,\cdots, \sigma_m$ of the bundle $E\to X$ such that for the
smooth sections
\begin{eqnarray*}
&&\Phi:X\times\R^m\to\Pi_1^\ast E,\;(y,{\bf t})\mapsto
S(y)+\sum^m_{i=1}t_i\sigma_i(y),\\
&&\Phi_{\bf t}:X\to E,\;y\mapsto S(y)+\sum^m_{i=1}t_i\sigma_i(y),
 \end{eqnarray*}
  where ${\bf t}=(t_1,\cdots, t_m)\in\R^m$ and $\Pi_1$ is the projection to the first factor
 of $X\times\R^m$, the following holds:
 There exist an open neighborhood ${\cal W}\subset {\cal O}(Z(S))$
of $Z(S)$ and  a small $\varepsilon>0$ such that:
\begin{description}
\item[(A)] The zero locus of $\Phi$ in $Cl({\cal W}\times
B_\varepsilon(\R^m))$ is compact. Consequently,  for any given small
open neighborhood ${\mathcal U}$ of $Z(S)$ there exists a
$\epsilon\in (0, \varepsilon]$ such that $Cl({\cal
W})\cap\Phi^{-1}_{{\bf t}}(0)\subset{\mathcal U}$ for any ${\bf
t}\in B_\epsilon(\R^m)$. In particular,  each set ${\cal
W}\cap\Phi^{-1}_{\bf t}(0)$ is {\rm compact} for ${\bf t}\in
B_\varepsilon(\R^m)$ sufficiently small.

\item[(B)] The restriction of $\Phi$ to ${\cal W}\times
B_\varepsilon(\R^m)$ is (strong) Fredholm and also transversal to
the zero section. So
$$ U_\varepsilon:=\{(y,{\bf t})\in {\cal W}\times
B_\varepsilon(\R^m)\,|\, \Phi(y,{\bf t})=0\}
$$
is a smooth manifold of dimension $m+{\rm Ind}(S)$, and  for
 ${\bf t}\in B_\varepsilon(\R^m)$ the section
 $\Phi_{\bf t}|_{\cal W}: X\to E$
 is transversal to the zero section if and only if ${\bf t}$ is
a regular value of the (proper) projection
 $$P_\varepsilon: U_\varepsilon\to B_\varepsilon(\R^m),\;(y,{\bf t})\mapsto{\bf t},$$
 and $\Phi_{\bf t}^{-1}(0)\cap{\cal W}=P_\varepsilon^{-1}({\bf t})$.
(Specially, ${\bf t}=0$ is a regular value of $P_\varepsilon$ if $S$
is transversal to the zero section).
 Then the Sard theorem yields  a residual subset
  $B_\varepsilon(\R^m)_{res}\subset B_\varepsilon(\R^m)$  such that:
 \begin{description}
 \item[(B.1)] For each
 ${\bf t}\in B_\varepsilon(\R^m)_{res}$ the set
  $(\Phi_{\bf t}|_{\cal W})^{-1}(0)\thickapprox(\Phi_{\bf t}|_{\cal W})^{-1}(0)\times\{{\bf t}\}
  =P_\varepsilon^{-1}({\bf t})$
  is a compact smooth manifold of dimension
  ${\rm Ind}(S)$ and  all $k$-boundaries
$$\partial^k(\Phi_{\bf t}|_{\cal W})^{-1}(0)=(\partial^k X)\cap (\Phi_{\bf t}|_{\cal W})^{-1}(0)$$
for $k=1, 2,\cdots$. Specially, if $Z(S)\subset {\rm Int}(X)$ one
can shrink $\varepsilon>0$ so that $(\Phi_{\bf t}|_{\cal
W})^{-1}(0)$ is a closed manifold for each ${\bf t}\in
B_\varepsilon(\R^m)_{res}$.

\item[(B.2)] If the Banach Fredholm bundle $(X,E, S)$ is {\bf
oriented}, i.e., the determinant  bundle ${\rm det}(DS)\to Z(S)$ is
given  a nowhere vanishing continuous section over $Z(S)$, then it
determines an orientation on $U_\varepsilon$. In particular, it
induces a natural orientation on every $(\Phi_{\bf t}|_{\cal
W})^{-1}(0)$ for ${\bf t}\in B_\varepsilon(\R^m)_{res}$.

 \item[(B.3)]   For any $l\in\N$ and two different  ${\bf t}^{(1)},
  {\bf t}^{(2)}\in B_\varepsilon(\R^m)_{res}$ the smooth
  manifolds $(\Phi_{{\bf t}^{(1)}}|_{\cal W})^{-1}(0)$ and $(\Phi_{{\bf
  t}^{(2)}}|_{\cal W})^{-1}(0)$ are cobordant in the sense that for a generic
  $C^l$-path $\gamma:[0,1]\to
B_\varepsilon(\R^m)$ with $\gamma(0)={\bf t}^{(1)}$ and
$\gamma(1)={\bf t}^{(2)}$ the set
$$\Phi^{-1}(\gamma):=\cup_{t\in
[0,1]}\{t\}\times (\Phi_{\gamma(t)}|_{\cal W})^{-1}(0)
$$
is a compact smooth manifold with boundary
$$\{0\}\times
(\Phi_{{\bf t}^{(1)}}|_{\cal W})^{-1}(0)\cup(-\{1\}\times(\Phi_{{\bf
t}^{(2)}}|_{\cal W})^{-1}(0)).$$
 In particular,  if $Z(S)\subset {\rm Int}(X)$ and
 $\varepsilon>0$ is suitably shrunk so that $(\Phi_{\bf t}|_{\cal W})^{-1}(0)\subset {\rm Int}(X)$
for any ${\bf t}\in B_\varepsilon(\R^m)$ then $\Phi^{-1}(\gamma)$
has no corners.

\item[(B.4)] The cobordant class of the manifold $(\Phi_{\bf
t}|_{\cal W})^{-1}(0)$ above is independent of all related choices.
  \end{description}
  \end{description}
  \end{theorem}

Now we furthermore assume that  $N$ is a connected manifold of
dimension $r$ and $f:X\to N$ is a smooth map.  When $X$ has no
boundary, by Theorem~\ref{th:2.9}(B.1), for each  ${\bf t}\in
B_\varepsilon(\R^m)_{res}$ the section $\Phi_{\bf t}:X\to E$ is
transversal to the zero section and the set
  $(\Phi_{\bf t}|_{\cal W})^{-1}(0)\subset X$
  is a compact smooth manifold of dimension $r$ and without
  boundary. So we may consider the $\Z_2$-Brouwer degree
$$
 \deg_{\Z_2}(f|_{(\Phi_{\bf t}|_{\cal W})^{-1}(0)})
 $$
of the restriction $f|_{(\Phi_{\bf t}|_{\cal W})^{-1}(0)}:
(\Phi_{\bf t}|_{\cal W})^{-1}(0)\to N$. The elementary properties
and Theorem~\ref{th:2.9}(B.3) show that $\deg_{\Z_2}(f|_{(\Phi_{\bf
t}|_{\cal W})^{-1}(0)})\in\Z_2$ is independent of the choice of
${\bf t}\in B_\varepsilon(\R^m)_{res}$. Moreover, it is claimed in
Theorem~\ref{th:2.9}(B.4) that the cobordant class of the manifold
$(\Phi_{\bf t}|_{\cal W})^{-1}(0)$ above is independent of all
related choices. Namely, suppose that  $\sigma'_1, \sigma'_2,\cdots,
\sigma'_{m'}$ are another group of smooth sections of the bundle
$E\to X$ such that the  section
$$
\Psi: {\cal W}'\times B_{\varepsilon'}(\R^{m'})\to \Pi^\ast_1E,\;
(y, {\bf t}')\mapsto S(y)+ \sum^{m'}_{i=1}t'_i\sigma'_i(y),
$$
 is Fredholm and transversal to the zero and that  the  set $\Psi^{-1}_{{\bf
t}'}(0)$ is compact for each ${\bf t}'\in
B_{\varepsilon'}(\R^{m'})$, where the section $\Psi_{{\bf t}'}:
{\cal W}'\to E$ is given by $\Psi_{{\bf t}'}(y)=\Psi(y,{\bf t}')$.
Let $B_{\varepsilon'}(\R^{m'})_{res}\subset
B_{\varepsilon'}(\R^{m'})$ be the corresponding residual subset such
that for each ${\bf t}'\in B_{\varepsilon'}(\R^{m'})_{res}$ the
section $\Psi_{{\bf t}'}$ is transversal to the zero section and
that any two ${\bf t}', {\bf s}'\in B_{\varepsilon'}(\R^{m'})_{res}$
yield cobordant manifolds $(\Psi_{{\bf t}'})^{-1}(0)$  and
$(\Psi_{{\bf s}'})^{-1}(0)$. Then it was shown in the proof of
\cite[Theorem 1.5(B.4)]{LuT} that there exist a compact submanifold
$\Theta_{({\bf t},{\bf t}')}^{-1}(0)\subset X\times [0,1]$ of
dimension $r+1$ for any ${\bf t}\in B^{reg}_{\varepsilon}(\R^m)$ and
${\bf t}'\in B^{reg}_{\varepsilon'}(\R^{m'})$ such that
$\partial\Theta_{({\bf t},{\bf t}')}^{-1}(0)=(\Phi_{\bf t}|_{\cal
W})^{-1}(0)\times\{0\}\cup\Psi^{-1}_{{\bf t}'}(0)\times\{1\}$. This
implies that
$$
 \deg_{\Z_2}(f|_{(\Phi_{\bf t}|_{\cal W})^{-1}(0)})=\deg_{\Z_2}(f|_{(\Psi_{{\bf t}'}|_{{\cal
 W}'})^{-1}(0)}).
 $$
Hence we have a well-defined $\Z_2$-value degree
\begin{equation}\label{e:3.1}
\deg_{\Z_2}(f, N, X, E, S):=\deg_{\Z_2}(f|_{(\Phi_{\bf t}|_{\cal
W})^{-1}(0)})\in\Z_2
\end{equation}
for any ${\bf t}\in B^{reg}_{\varepsilon}(\R^m)$, and call it {\bf
$\Z_2$-degree of $f:X\to N$ relative to $(X, E, S)$.} Of course,
when both $(X, E, S)$ and $N$ are oriented, we may define {\bf $\Z$-
degree of $f:X\to N$ relative to $(X, E, S)$.}

Let $\{S_\lambda\}_{\lambda\in [0, 1]}$ be a smooth family of smooth
Fredholm sections of the bundle $E\to X$ of index $r$ and with
compact zero sets. Then we can still choose finitely many smooth
sections $\sigma_1, \sigma_2,\cdots, \sigma_m$ of the bundle $E\to
X$, an open neighborhood ${\cal W}_\lambda$ of each
$Z(S_\lambda)\subset X$, and a residual subset
$B_{\varepsilon}(\R^{m})_{res}$ for some small $\varepsilon>0$, such
that for each ${\bf t}\in B_{\varepsilon}(\R^{m})_{res}$ the
restrictions of the smooth sections
\begin{eqnarray*}
&&\Phi^0_{\bf t}:X\to E,\;y\mapsto S_0(y)+\sum^m_{i=1}t_i\sigma_i(y),\\
&&\Phi^1_{\bf t}:X\to E,\;y\mapsto
S_1(y)+\sum^m_{i=1}t_i\sigma_i(y),\\
&&\Phi_{\bf t}:X\times [0, 1]\to\Pi_1^\ast E,\;(y, \lambda)\mapsto
S_\lambda(y)+\sum^m_{i=1}t_i\sigma_i(y)
 \end{eqnarray*}
to  ${\cal W}_0$, ${\cal W}_1$ and ${\cal W}=\cup_{\lambda\in
[0,1]}{\cal W}_\lambda$ are transversal to the zero sections
respectively. In particular, we get
$$
\partial (\Phi_{\bf t}|_{{\cal
W}})^{-1}(0)=(\Phi^0_{\bf t}|_{{\cal
W}_0})^{-1}(0)\times\{0\}\bigcup(\Phi^1_{{\bf t}}|_{{\cal
W}_1})^{-1}(0)\times\{1\}.
$$
It follows that
\begin{equation}\label{e:3.2}
\deg_{\Z_2}(f, N, X, E, S_0)=\deg_{\Z_2}(f, N, X, E, S_1)
\end{equation}
and thus $\deg_{\Z_2}(f, N, X, E, S_\lambda)$ {\bf is independent
of} $\lambda\in [0,1]$.

Similarly,  if $(X, E, S_\lambda)$ and $N$ are oriented,
$\deg_{\Z}(f, N, X, E, S_\lambda)$ is independent of $\lambda\in
[0,1]$ as well. \vspace{2mm}

\noindent{\bf Proof of Theorem~\ref{th:2.7}.}\quad Define the
evaluation map
\begin{equation}\label{e:3.3}
\Theta:{\cal
B}^\tau\to L,\;u\mapsto u(1),
\end{equation}
where
$1\in \C\subset\C\cup\{\infty\}=S^2$. Applying the arguments above
to  the Banach Fredholm bundle $({\cal B}^\tau, {\cal E}^+, {\cal
F}_{T,\lambda}|_{{\cal B}^\tau})$, $\lambda\in [0,1]$, we arrive at
\begin{equation}\label{e:3.4}
\deg_{\Z_2}(\Theta, L, {\cal B}^\tau, {\cal E}^+, {\cal
F}_{T,1}|_{{\cal B}^\tau})=\deg_{\Z_2}(\Theta, L, {\cal B}^\tau,
{\cal E}^+, {\cal F}_{T,0}|_{{\cal B}^\tau})
\end{equation}
by (\ref{e:3.2}).
 Since each $w\in{\cal B}$ is contractible, ${\cal Z}^\tau_{T,0}=({\cal
F}_{T,0}|_{{\cal B}^\tau})^{-1}(0_{{\cal E}^+})$ precisely consists
of the constant maps $S^2\to L$. It is easily proved that ${\cal
F}_{T,0}|_{{\cal B}^\tau}:{\cal B}^\tau\to {\cal E}^+$ is
transversal to the zero section, and that (\ref{e:3.1}) yields
\begin{equation}\label{e:3.5}
\deg_{\Z_2}(\Theta, L, {\cal B}^\tau, {\cal E}^+, {\cal
F}_{T,0}|_{{\cal B}^\tau})=1.
\end{equation}
Let $F$ be a smooth perturbation section  of ${\cal F}_{T,1}|_{{\cal
B}^\tau}$ as $\Phi^1_{\bf t}$ above. Choose $l_0\in L$ to be a
regular value for the evaluations
$$
\Theta|_{F^{-1}(0_{{\cal E}^+})}: F^{-1}(0_{{\cal E}^+})\to L.
$$
Then (\ref{e:3.4}) and (\ref{e:3.5}) show that
$$
 \deg_{\Z_2}(\Theta|_{F^{-1}(0_{{\cal E}^+})}, l_0)=1.
$$
Hence $\Theta|_{F^{-1}(0_{{\cal E}^+})}: F^{-1}(0_{{\cal E}^+})\to
L$ induces an injection map
\begin{equation}\label{e:3.6}
(\Theta|_{F^{-1}(0_{{\cal E}^+})})^\ast: \check{H}^\ast(L,\Z_2)\to
\check{H}^\ast(F^{-1}(0_{{\cal E}^+}),\Z_2).
\end{equation}

Note that  $F^{-1}(0_{{\cal E}^+})$ can be chosen so close to ${\cal
Z}^\tau_{T,1}$ that it is contained a given small neighborhood of
${\cal Z}^\tau_{T,1}$ for which Theorem~\ref{th:2.6} implies for
$T\ge T_0>6$
\begin{equation}\label{e:3.7}
\sigma_T(u|_{Z_T})\in U\quad\forall w\in F^{-1}(0_{{\cal
E}^+})\;{\rm and}\;u=w\circ\phi.
\end{equation}
  Here we use $F^{-1}(0_{{\cal E}^+})\subset C_c^\infty(S^2, M)$ due to the arguments above
 Lemma~\ref{lem:2.1}.
 Define
$$
\Xi:F^{-1}(0_{{\cal E}^+})\to X^\tau_{T,d},\; w\mapsto
u|_{Z_T}\;{\rm for}\;u=w\circ\phi,
$$
 by (\ref{e:2.22}), (\ref{e:2.24}), (\ref{e:3.3}) and (\ref{e:3.7})
 it is easy to see that we have for $T\ge T_0$ the commutative
diagram
\begin{center}\setlength{\unitlength}{1mm}
\begin{picture}(80,30)
\thinlines \put(10,25){$X^\tau_{T}$} \put(22,25){\vector(1,0){43}}
\put(35,27){$\sigma_T$} \put(15,6){\vector(0,1){16}}
\put(10,13){$\Xi$} \put(70,25){$U$} \put(70,22){\vector(0,-1){16}}
\put(63,13){$\pi|_U$} \put(2 ,0){$F^{-1}(0_{{\cal E}^+})$}
\put(22,1){\vector(1,0){43}} \put(35,3){$\Theta|_{F^{-1}(0_{{\cal
E}^+})}$} \put(70,0){$L$}
\end{picture}\end{center}
By (\ref{e:3.6}) we get the injectiveness of the map
$$
(\pi|_U)^\ast: \check{H}^\ast(L,\Z_2)\to \check{H}^\ast(U,\Z_2).
$$

If $L$ is orientable,  the Banach Fredholm bundles $({\cal B}^\tau,
{\cal E}^+, {\cal F}_{T, 0}|_{{\cal B}^\tau})$, and therefore
$({\cal B}^\tau, {\cal E}^+, {\cal F}_{T,\lambda}|_{{\cal
B}^\tau})$, $\lambda\in [0,1]$, are orientable. In this case we can
define $\Z$-degree $\deg_{\Z}(\Theta, L, {\cal B}^\tau, {\cal E}^+,
{\cal F}_{T,\lambda}|_{{\cal B}^\tau})$ and get $\deg_{\Z}(\Theta,
L, {\cal B}^\tau, {\cal E}^+, {\cal F}_{T,\lambda}|_{{\cal
B}^\tau})\in\{1, -1\}$. The desired conclusion follows immediately.
\hfill$\Box$\vspace{2mm}

\section{Examples and  further programme}\label{sec:4}
\setcounter{equation}{0}

\begin{example}\label{ex:4.1}
{\rm Consider the torus $T^{2n}=\R^{2n}/\Z^{2n}$ with the standard symplectic form $\omega=dx\wedge dy)$.
There exist three natural anti-symplectic involutions on it given by\\
$\bullet$ $\tau_1:T^{2n}\to T^{2n},\;[x,y]\mapsto [-x,y]$,\\
$\bullet$ $\tau_2:T^{2n}\to T^{2n},\;[x,y]\mapsto [x,-y]$,\\
$\bullet$ $\tau_3:T^{2n}\to T^{2n},\;[x,y]\mapsto [y, x]$.\\
Clearly, ${\rm Fix}(\tau_1)=[0]\times T^n\subset T^{2n}$, ${\rm Fix}(\tau_2)=T^n\times [0]\subset T^{2n}$
and ${\rm Fix}(\tau_3)=\{[x,x]\in T^{2n}\,|\,x\in\mathbb{R}^n\}$.
 Let $H\in C^\infty(\R\times\R^{2n},\R)$ be $1$-periodic in
all its variables so that it may be viewed as a Hamiltonian function
on the standard torus which is 1-periodic in time. Corresponding to the three cases above, suppose that Hamiltonian $H$ also, respectively, satisfies\\
$\bullet$ $H(-t, x,y)=H(t, -x, y)$ for any $t\in\R$ and $z=(x,y)\in\R^{2n}$,\\
$\bullet$ $H(-t, x,y)=H(t, x, -y)$) for any $t\in\R$ and $z=(x,y)\in\R^{2n}$,\\
$\bullet$ $H(-t,x,y)=H(t, y,x)$ for any $t\in\R$ and $z=(x,y)\in\R^{2n}$.\\
Then in each case Theorem~\ref{th:1.3} gives at least $n+1$ contractible $1$-periodic solutions
$\gamma:\R\to T^{2n}$ of $\dot\gamma(t)=X_{H}(t,\gamma(t))$ satisfying
$\gamma(-t)=\tau_i(\gamma(t))$ for any $t\in\R$, $i=1,2,3$, respectively.
It is the contractibility of $\gamma$ that there exists a lift loop
$z=(x,y):\R/\Z\to\R^{2n}$ of it satisfying
 the associated Hamiltonian system on $\mathbb{R}^{2n}$
$$
\dot z=J\nabla H(t, z)\quad \hbox{with}\;J=\left(\begin{array}{cr}
0& -I_n\\
I_n& 0\end{array}\right)
$$
and the following conditions\\
$\bullet$  $x(-t)+x(t)\in\mathbb{Z}^n$ and $y(-t)-y(t)\in\mathbb{Z}^n$ for any $t\in\R$,\\
$\bullet$  $x(-t)-x(t)\in\mathbb{Z}^n$ and $y(-t)+y(t)\in\mathbb{Z}^n$ for any $t\in\R$,\\
$\bullet$  $x(-t)-y(t)\in\mathbb{Z}^n$ and $y(-t)-x(t)\in\mathbb{Z}^n$ for any $t\in\R$,\\
respectively.}
\end{example}

Clearly, the conclusions in Example~\ref{ex:4.1}  cannot be derived from  \cite[Th.1]{CoZe} though the latter yields at $2n+1$ periodic solutions of
$\dot z=J\nabla H(t, z)$ of period $1$.

\begin{example}\label{ex:4.2}{\rm
{\bf (i)} Consider the standard complex projective space $\CP^{n}$
with the Fubini-Study form $\omega_{FS}$ satisfying $\int_{\CP^{1}}\omega_{\rm FS}=\pi$.
Then the rationality index of $(\CP^{n}, \omega_{\rm FS})$ is equal to $\pi$.
 Let $H\in C^\infty(\R\times\CP^{n},\R)$ be $1$-periodic
in the first variable, and also satisfy $H(-t, [z])=H(t,
\sigma([z]))$ for any $t\in\R$ and $[z]\in\CP^n$, where $\sigma$ is
the standard complex conjugation on $\CP^n$ with ${\rm Fix}(\sigma)=\RP^n$.
Then the associated Hamiltonian system $\dot z=X_H(t, z)$ on
$(\CP^n,\omega_{\rm FS})$ has at least $n+1$ contractible periodic solutions
$z:\R\to\CP^n$ of period $1$ satisfying $z(-t)=\sigma(z(t))$ for any
$t\in\R$ provided the Hofer norm $\|H\|<\pi$ by Theorem~\ref{th:1.3}.
(Actually, the final restriction ``$\|H\|<\pi$''
may be removed out with Fortune's method in \cite{Fo}.)\\
{\bf (ii)} Let $(P, \omega)$ be a simply connected closed symplectic
manifold of dimension $4$ and with $c_1(TP)|_{\pi_2(P)}=0$.  By the
Hurewicz isomorphism theorem and the Poincar\'e dual theorem there
exists a class $A\in\pi_2(P)$ such that $\beta(A)>0$. So
$r_\omega\in [0, +\infty)$. It easily follows from
\cite[Theorem 3.1.5]{McSa2} that for generic $J\in{\cal J}(P,\omega)$
there is no nonconstant $J$-holomorphic spheres in $P$, and thus
$m(P,\omega, J)=+\infty$. Moreover, if $(P, \omega)$ is also
real symplectic,  by \cite[Proposition
11.10]{FuOOO} for generic $J\in\R{\cal J}(P,\omega)$ there is no nonconstant $J$-holomorphic sphere in $P$
 and so $m(P,\omega, J)=+\infty$. A well-known
example of such real symplectic manifolds is
the $K3$-surface
$$
X=\bigl\{[z_0:\cdots:z_3]\in\CP^3\,\bigm|\,
\sum^3_{j=0}z_j^4=0\bigr\}
$$
with the canonical symplectic structure $\omega_{\rm can}$ induced by
the form $\omega_{FS}$ on $\CP^n$ as in (i) above and
with the anti-symplectic involution induced by  the
standard complex conjugation on $\CP^n$ (cf. \cite[Example 4.27]{McSa1}).
Hence $\pi\le r_{\omega_{\rm can}}<+\infty$. Note that the real part of $X$ is empty!\\
{\bf (iii)} A symplectic manifold $(M, \omega)$ of dimension $2n$ is
said to be {\it negative monotone} if
$c_1(TM)|_{\pi_2(M)}=\lambda\cdot\omega|_{\pi_2(M)}$ for some
negative constant $\lambda$, and {\it semipositive} if either
$\omega(M)|_{\pi_2(M)}=\mu\cdot c_1|_{\pi_2(M)}$ for some constant
$\mu\ge0$, or $ c_1|_{\pi_2(M)}=0$ or the minimal Chern number $N\ge
n-2$, see \cite[Exercise 6.4.3]{McSa2}. Here the minimal Chern
number $N$ of $(M, \omega)$ is the positive generator of
$c_1(M)\bigl(\pi_2(M)\bigr)$ if $c_1|_{\pi_2(M)}\ne 0$, and
$+\infty$ if $c_1|_{\pi_2(M)}=0$.
Note that a simply connected and closed  symplectic manifold
has always finite rationality index by
the Hurewicz isomorphism theorem and the universal coefficient theorem.
In a negative monotone symplectic manifold $(M, \omega)$ with minimal
Chern number $N\ge \dim M/2$, for generic $J\in{\cal J}(M, \omega)$ there is no
nonconstant $J$-holomorphic sphere  and hence $m(M, \omega,
J)=+\infty$ by \cite{McSa2}, and  $m(M, \omega, J)=+\infty$ for
generic $J\in\R{\cal J}(M, \omega)$ by  \cite[Proposition
11.10]{FuOOO}  if $(M, \omega)$ is also real.   {\bf Here are some
concrete examples}, which were in details discussed in
\cite[Appendix A]{Laz}. For an integer $n\ge 4$ and an odd integer $d$
let
$$
M_{n,d}=\bigl\{[z_0:\cdots:z_n]\in\CP^n\,\bigm|\,
\sum^3_{j=0}z_j^d=0\bigr\}
$$
equipped with a canonical symplectic structure $\omega_{n,d}$ induced by
the form $\omega_{FS}$ on $\CP^n$ as in (i) above.
Then $\pi\le r_{\omega_{n,d}}<+\infty$. The standard complex conjugation on $\CP^n$
induces an anti-symplectic $\tau$ on $M_{n,d}$ with ${\rm Fix}(\tau)=M_{n,d}\cap\RP^n$ which is
homeomorphic to $\RP^{n-1}$. So if $H\in C^\infty(\R\times M_{n,d},\R)$ is $1$-periodic
in the first variable, and also satisfy $H(-t, [z])=H(t,
\tau([z]))$ for any $t\in\R$ and $[z]\in M_{n,d}$, we have
$\sharp {\cal P}_0(H,\tau)\ge n$ provided $\|H\|<\pi$.

 It was shown in \cite[Appendix A]{Laz} that $M_{n,d}$ is simply connected, has a
minimal Chern number $N_{n,d}=|n+1-d|$, and satisfies
$$
c_1(M_{n,d})|_{\pi_2(M_{n,d})}=\frac{n+1-d}{r}\cdot\omega_{n,d}|_{\pi_2(M_{n,d})}
$$
for some $r>0$. Since $\dim M_{n,d}=2n-2$, $M_{n,d}$ is negative
monotone if and only if $n+1<d$, and $N_{n,d}\ge\frac{1}{2}\dim
M_{n,d}$ if and only if $d\ge 2n$ or $d=1$. Hence the arguments above
show that each $M_{n,d}$ with $n\ge 4$ and odd integer $d\ge 2n$ or $d=1$
satisfies
$$
m(M_{n,d},\omega_{n,d}, J)=+\infty
$$
for generic $J\in{\cal J}(M_{n,d},\omega_{n,d})$ and for generic $J\in\R{\cal J}(M_{n,d},\omega_{n,d})$.   Nonsingular algebraic subvarieties of $\CP^n$ defined by real equations
may provide more examples.}
\end{example}

{\bf Our programme} \cite{Lu3} is to construct a real Floer homology
$FH_\ast(M, \omega, \tau, H)$ for a real symplectic manifold $(M,
\omega, \tau)$ with nonempty compact $L={\rm Fix}(\tau)$ only using
${\cal P}_0(H, \tau)$, which may be viewed as an {\it intermediate}
between the Floer homology for Hamiltionian maps and the Floer
homology for Lagrangian intersections,  to prove that it is
isomorphic to $H_\ast(L)\otimes R_\omega$ for some Novikov ring
$R_\omega$, and then to relate it to some possible open
GW-invariants and something
 as in \cite{FuOOO}, \cite{BiCo} and Auroux's talk at Montreal, May 19-24, 2008.

%\appendix
%\section{ A}\label{app:A}\setcounter{equation}{0}

\end{document}